\documentclass[12pt]{article}
\usepackage{amsfonts}
\usepackage{amsmath}
\usepackage{amssymb}
\usepackage{amsthm}
\usepackage{epsfig}
\theoremstyle{plain}
\newtheorem{theo}{Theorem}[section]
\newtheorem{pro}[theo]{Proposition}
\newtheorem{lem}[theo]{Lemma}
\newtheorem{cor}[theo]{Corollary}
\theoremstyle{definition}
\newtheorem{defn}[theo]{Definition}
\theoremstyle{remark}
\newtheorem{rem}[theo]{Remark}
\newenvironment{demo}[1]{\par\begin{trivlist}%
\item[]{\bf #1}\ }{\end{trivlist}\par}
\newcommand{\Proof}{\begin{demo}{{\it Proof.\ }}}
\newcommand{\Proofwithname}[1]{\begin{demo}{{\it #1}}}
\newcommand{\QED}{\qed\end{demo}}

\newcommand{\bl}{L}
\newcommand{\ca}{{\mathcal A}}
\newcommand{\cb}{{\mathcal B}}
\newcommand{\cc}{{\mathcal C}}

\newcommand{\ce}{{\mathcal E}}
\newcommand{\cf}{{\mathcal F}}
\newcommand{\cg}{{\mathcal G}}
\newcommand{\ch}{{\mathcal H}}
\newcommand{\ck}{{\mathcal K}}

\newcommand{\cn}{{\mathcal N}}
\newcommand{\bz}{{\mathbb Z}}
\newcommand{\br}{{\mathbb R}}
\newcommand{\bn}{{\mathbb N}}
\newcommand{\bc}{{\mathbb C}}

\newcommand{\fg}{{\mathfrak G}}
\newcommand{\tce}{\tilde {\cal E}}
\newcommand{\tcf}{\tilde {\cal F}}
\newcommand{\tmu}{\tilde{\mu}}
\newcommand{\tg}{\tilde{G}}
\newcommand{\tx}{\tilde{X}}
\newcommand{\lr}{\leftrightarrow}
\newcommand{\esssup}{\mathop{\rm esssup}}
\newcommand{\essinf}{\mathop{\rm essinf}}
\newcommand{\supp}{\mathop{\rm supp}}
\newcommand{\Cp}{\mathop{\rm Cap}}
\newcommand{\vp}{\varphi}
\newcommand{\sd}{\mathsf d}

\def\mint{{}-\!\!\!\!\!\!\!\int}
\def\dl{\delta}
\def\lll{\|}
\def\qed{\hbox {\rlap{$\sqcup$}{$\sqcap$}}}
\makeatletter
\def\rom#1{\mbox{\leavevmode\skip@\lastskip\unskip\/%
\ifdim\skip@=\z@\else\hskip\skip@\fi{\rm{#1}}}}
\@addtoreset{equation}{section}

\makeatother
\begin{document}
\title{A trace theorem for Dirichlet forms on fractals}
\author{Masanori Hino%
\thanks{Research partially supported by Ministry of Education,
Culture, Sports, Science and Technology, Grant-in-Aid for Encouragement 
of Young Scientists, 15740089.}\\
Graduate School of Informatics,\\
Kyoto University,\\ 
Kyoto 606-8501, Japan\\
{\tt hino@i.kyoto-u.ac.jp}
\and
Takashi Kumagai%
\thanks{Research partially supported by Ministry of Education,
Culture, Sports, Science and Technology, Grant-in-Aid for Encouragement 
of Young Scientists, 16740052.}
\footnote{Corresponding author. Tel: +81-75-753-7245~ 
Fax: +81-75-753-7272}
\\Research Institute for Mathematical\\ Sciences,  
Kyoto University,\\ Kyoto 606-8502, Japan\\
{\tt kumagai@kurims.kyoto-u.ac.jp}
\and\\
\small{Dedicated to Professor S. Watanabe on the occasion of his 
70th birthday.}}

\date{19 October, 2005}

\maketitle
\begin{abstract}
We consider a trace theorem for self-similar Dirichlet forms on 
self-similar sets to self-similar subsets. In particular, we characterize 
the trace of the domains of Dirichlet forms on the Sierpinski gaskets and 
the Sierpinski carpets to their boundaries, where boundaries mean the triangles 
and rectangles which confine gaskets and carpets. As an application, we 
construct diffusion processes on a collection of fractals called fractal fields,
which behave as the appropriate fractal diffusion within each fractal component
of the field.  
\smallskip

\noindent{\em MSC:} 46E35; 28A80; 31C25; 60J60 

\noindent{\em Keywords:} Trace theorem; Self-similar sets; Dirichlet forms; 
Diffusions on fractals; Lipschitz spaces; Besov spaces; Sierpinski carpets 

\end{abstract}

\newpage
\section{Introduction}
The trace of Sobolev spaces on $\br^n$ to linear subspaces have been studied 
in various directions as generalizations of the Sobolev imbedding theorem. 
There has also been extensive study how to extend Sobolev, Besov and Lipschitz spaces
from subdomains of $\br^n$ to the whole spaces (see for example, \cite{adh,ste} and
the references therein). 
Since 80's, there are generalizations of these problems 
for Besov-type spaces on more complicated spaces, namely on the so-called Alfors
$d$-regular sets (\cite{jw,tr}). 

On the other hand, recent developments of analysis on fractals give new lights 
to these problems. On many fractals such as Sierpinski gaskets and Sierpinski carpets,
diffusion processes and the ``Laplace'' operators are constructed.
It turns out that the domains of 
the corresponding Dirichlet forms are Besov-Lipschitz spaces. 

\begin{figure}[ht]
\centerline{\epsfig{file=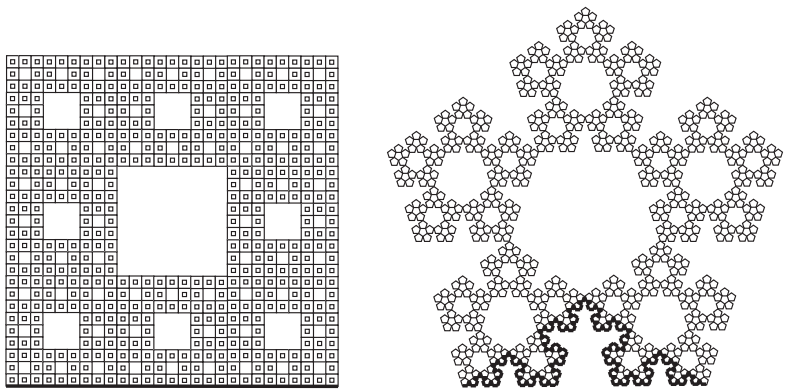, height=2in}}
\caption{The Sierpinski carpet and the Pentakun}
\end{figure}

In this paper, we consider the following natural question: given a Besov-type space 
on a self-similar fractal $K$, what is the trace of the space to a self-similar 
subspace $L$? We would indicate two examples in Figure 1. The left figure is 
when $K$ is the so-called 2-dimensional Sierpinski carpet (see Section~5~3) for
the definition) and $L$ is the line on the bottom (drawn by the thick line). The right 
figure is when $K$ is the Pentakun (a self-similar fractal determined by five contraction 
maps; see Section~5~2) for the definition) and $L$ is a Koch-like curve (drawn by the thick 
curve). In each case, the domain of the Dirichlet form on $K$ is the Besov-Lipschitz 
space, but one cannot obtain the trace using the general theory given by Jonsson-Wallin 
(\cite{jw}) and Triebel (\cite{tr}). 

This problem was quite recently solved by Jonsson (\cite{jo}) for one typical
case, i.e. when $K$ is the 2-dimensional Sierpinski gasket and $L$ is the bottom line. But his methods
rely strongly on the structure of the Sierpinski gasket and its Dirichlet form, and 
they cannot be applied to the so-called infinitely ramified fractals such as Sierpinski
carpets. Instead, we use the self-similarity of the form and some kind of uniform property 
of harmonic functions which can be guaranteed by the Harnack inequalities. 
Our methods can be applied to the Sierpinski carpets (even to the high dimensional ones) and
we can state the trace theorem under some abstract framework. In fact, we would need
various assumptions for $K$ and for the Dirichlet form on $K$, which are stated in Section~2. 
Unless these conditions are satisfied, there may be various possibility of the trace, because of the 
\lq\lq complexity\rq\rq\ 
of the space (see Section~5~4) for an example). 

In order to prove our trace theorem, we give a discrete
approximation of our Besov-Lipschitz space in Section~3.1. This approximation result is also new
and is regarded as a generalization of
the main result in \cite{kam}. 
The restriction theorem is given in Section~3.2; the key estimate   
(Proposition \ref{theo:keyhin}) is based on the idea used by one of the author in \cite{hin}.
The extension theorem is given in Section~3.3, where the classical 
construction 
of the Whitney decomposition and the extension map is 
modified and 
generalized to this framework.
 
\begin{figure}[ht]
\centerline{\epsfig{file=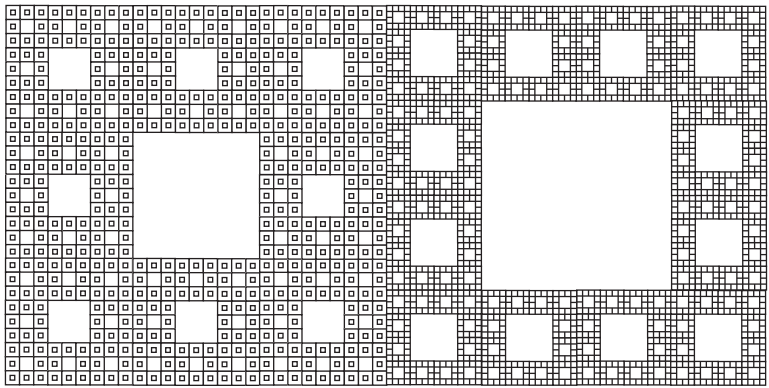, height=2in}}
\caption{An example of fractal fields}
\end{figure}

Such a trace theorem has an important application to the penetrating process, which is discussed in 
Section~6. Let us indicate one concrete example. Given two types of Sierpinski carpets as in Figure 2
(the left carpet is determined by contraction maps with the contraction rate $1/3$ and there is one hole 
in the middle, while the right carpet is determined by contraction maps with the contraction rate $1/4$ 
and there is one bigger hole in the middle). On each carpet, one can construct a self-similar
diffusion; the question is whether one can construct a diffusion which 
behaves as the appropriate fractal diffusions within each carpet and which penetrates each fractal.
In order to construct such a diffusion by the superposition of Dirichlet forms on each carpet,  
the key problem is whether there is enough functions whose restriction to each carpet is in the domain
of each Dirichlet form. To answer this question, it is crucial to get 
the information of the trace of the Dirichlet form on 
each carpet to the line, which is the intersection of the two carpets. Indeed, 
when one of the author studied this problem on fractals in \cite{kum,hk}, he needed
a very strong assumption on each fractal because of the lack of the information of the trace. 
Our trace theorem can be applied here and we can construct penetrating processes on much wider class 
of fractals.

Throughout this article, if $f$ and $g$ depend on a variable $x$ ranging 
in a set $A$, $f\asymp g$ means that there exists
$C>0$ such that
$C^{-1}f(x)\le g(x)\le C\,f(x)$ for all $x\in A$.
We will use $c$, with or without
subscripts, to denote strictly
positive constants whose values are insignificant.
\section{Framework and the main theorem}
Let $(X,\sd)$ be a complete separable 
metric space. 
For $\alpha>1$ and a finite index set $W$,
let $\{F_i\}_{i\in W}$ be a family of {\it
$\alpha$-similitudes} on $X$, i.e.\ $\sd(F_i(x),F_i(y))=\alpha^{-1}\sd(x,y)$
for all $x,y\in X$. 
Let $S$ be a subset of $W$ and let $N$ denote the cardinality of $S$.
Since 
$\{F_i\}_{i\in S}$ is a family of contraction maps,
there exists a unique non-void compact set $K$ such that ${K}
=\bigcup_{i\in S}F_i ({K})$. 
We assume that $K$ is connected.
Note that $W$ will be needed in general when we define a self-similar subset $L$ below.
In various important examples such as 1), 3) in Section~5, we can take $W=S$. 

We will make the relation to the shift space. 
The one-sided shift space $\Sigma$ is defined by $\Sigma = W^{\bn}$. For $w \in \Sigma$, 
we denote the $i$-th element in the sequence by $w_i$ and write $w =w_1 w_2 w_3 \cdots$.
When $w\in W^n$, $|w|$ denotes $n$.
For $v\in W^m$ and $w\in W^n$, we define $v\cdot w\in W^{m+n}$ by $v\cdot 
w=v_1 v_2\cdots v_m w_1 w_2\cdots w_n$. For $A\subset W^m$ and $B\subset W^n$, $A\cdot B$ denotes $\{v\cdot w: v\in A,\ w\in B\}$.
The set $w\cdot A$ is defined as $\{w\}\cdot A$.
By definition, $W^0=\{\emptyset\}$ and $\emptyset\cdot A=A$.

Let $\fg$ be a group consisting of isometries on $K$.
We assume the following.
\begin{itemize}
\item For each $i\in W$, there exist $j=j(i)\in S$ and $\Psi_i\in\fg$ such that $F_i=F_j\circ \Psi_i$.
\item For each $(\Psi,\alpha)\in\fg\times S$, there exists $(\hat\Psi,\hat\alpha)\in\fg\times S$ such that
$
\Psi\circ F_\alpha = F_{\hat\alpha}\circ \hat\Psi
$.
\end{itemize}
Note that, when $W=S$, we can always take as $\fg$ the trivial group consisting of one element.
We write $F_{w_1\cdots w_n}=F_{w_1}\circ F_{w_2}\circ \cdots \circ F_{w_n}$
for $w=w_1w_2\cdots w_n$. 
We regard $F_\emptyset$ as an identity map.
For $w\in W^n$ and $A\subset W^n$ for some $n\in\bz_+$, define $K_w=F_w(K)$ and $K_A=\bigcup_{v\in A}K_v$.
\begin{lem}\label{lem:Phi}
There exist maps $\Phi:\bigcup_{n\in\bz_+}W^n \to \bigcup_{n\in\bz_+}S^n$ and $\Psi:\bigcup_{n\in\bz_+}W^n \to \fg$ such that $F_w=F_{\Phi(w)}\circ\Psi(w)$ for each $w\in \bigcup_{n\in\bz_+}W_n$.
In particular, $K_w=K_{\Phi(w)}$.
\end{lem}
\Proof
  Set $\Phi(\emptyset)=\emptyset$ and $\Psi(\emptyset)={}$the unit element of $\fg$.
  When $i\in W^1$, it suffices to set $\Phi(i)=j(i)$ and $\Psi(i)=\Psi_i$.
  Suppose that $\Phi(w)$ is defined for $w\in W^{n}$.
  Then, for $w'=w\cdot i$ with $i\in W$,
$
  F_{w'}=F_w\circ F_i=F_{\Phi(w)}\circ\Psi(w)\circ F_{j(i)}\circ\Psi_{i}
$.
  This is equal to $F_{\Phi(w)}\circ F_{\hat i}\circ\hat\Psi\circ \Psi_i$ for some $(\hat\Psi,\hat i)\in\fg\times
S$.
  Therefore, it is enough to define $\Phi(w')=\Phi(w)\cdot\hat i$ and $\Psi(w')=\hat\Psi\circ \Psi_i$.
\QED
Define $\pi:\Sigma\to K$ by the relation $\{\pi (w)\}=\bigcap_m K_{w_1\cdots w_m}$ for $w=w_1w_2\cdots\in\Sigma$. 
Define 
\begin{equation}\label{eq:critpc}
C_K := \pi^{-1} \left( \bigcup_{i,j \in S, i \ne j} (K_i \cap K_j)\right),\qquad
P_K := \bigcup_{n \ge 1} \sigma^n (C_K),\end{equation}
where $\sigma : \Sigma \to \Sigma$  is the left shift map, i.e.
$\sigma w = w_2 w_3 \cdots$ if $w = w_1 w_2 w_3\cdots$. 

For $v,w\in W^n$, we write $v\stackrel{n,K}{\sim} w$ if $K_v\cap K_w\ne \emptyset$.
For $w\in W^n$ and $A\subset W^n$, $w\stackrel{n,K}{\sim} A$ means that 
$w\stackrel{n,K}{\sim} v$ for some $v\in A$.
For $A\subset W^n$, define $\cn_0(A)=A$ and $\cn_k(A)=\{v\in W^n\mid 
v\stackrel{n,K}{\sim} \cn_{k-1}(A)\}$ for $k\in\bn$ inductively.
We set $\cn_k(w)=\cn_k(\{w\})$ for $w\in W^n$.

Let $I$ be a subset of $W$. 
We assume that the cardinality $N_I$ of $I$ is less than $N$.
Let $L$ be a unique non-void compact set such that $L
=\bigcup_{i\in I}F_i (L)$. 
Clearly, $L$ is a subset of $K$.
Denote $F_w(L)$ by $L_w$ for $w\in\bigcup_{n\in\bz_+}I^n$.
Let $M\in\bn$.
For $v$, $w\in I^n$, we write $v\underset{M}{\stackrel{n,L}{\longleftrightarrow}} w$ if $v \in \cn_M(w)$. 
We fix $M$ so that for each $i$, $j\in I$, there exist $i_1,i_2,\ldots\in I$ satisfying 
$i\underset{M}{\stackrel{1,L}{\longleftrightarrow}}i_1\underset{M}{\stackrel{1,L}
{\longleftrightarrow}}i_2\underset{M}{\stackrel{1,L}{\longleftrightarrow}}\cdots\underset{M}{\stackrel{1,L}{\longleftrightarrow}}j$.
In what follows, we omit $M$ from the notation $\underset{M}{\stackrel{n,L}{\longleftrightarrow}}$.
We assume the following.

\begin{enumerate}
\item[(A1)] $\sup_{n\in\bz_+}\max_{w\in S^n}\#(\cn_1(w)\cap S^n)<\infty$
and $C_0:=\sup_{n\in\bz_+}\max_{w\in I^n}\#(\cn_M(w)\cap I^n)<\infty$.
\item[(A2)] There exist $k_1, k_2>0$ such that, for 
$x,y\in L$, $n\in\bz_+$ and $v,w\in I^n$ with $x\in K_v$ and $y\in K_w$, $\sd(x,y)< k_1 \alpha^{-n}$ implies 
$v\stackrel{n,L}{\lr}w$ and $v\stackrel{n,L}{\lr}w$ implies $\sd(x,y)< k_2 \alpha^{-n}$.
\item[(A3)] There exist $k_1, k_2>0$ such that, for 
$x,y\in K$, $n\in\bz_+$ and $v,w\in S^n$ with $x\in K_v$ and $y\in K_w$, $\sd(x,y)< k_1 \alpha^{-n}$ implies 
$v\stackrel{n,K}{\sim}w$ and $v\stackrel{n,K}{\sim}w$ implies $\sd(x,y)< k_2 \alpha^{-n}$.
\end{enumerate}
Let $\hat\mu$ and $\hat\nu$ be the canonical Bernoulli measures on $S^\bn$ and $I^\bn$, respectively.
That is, they are infinite product measures of $S$ (resp.\ $I$) with uniformly distributed measure.
Denote by $\mu$ the image measures of $\hat\mu$ by the map $\pi|_{S^\bn}:S^\bn\to K$.
In the same way, the probability measure $\nu$ on $L$ is defined.
By conditions (A1), (A2), and (A3) and \cite[Theorem~1.5.7]{kig}, the Hausdorff dimensions of $K$ and $L$ are equal to $d_f:=\log N/\log\alpha$ and $d:=\log N_I/\log\alpha$, respectively, and $\mu$ and $\nu$ are equivalent to the Hausdorff measures on $K$ and $L$, respectively.

We will further assume the following.
\begin{enumerate}
\item[(A4)] $\mu(\{x\in K: \#
(\pi^{-1}(x)\cap S^\bn)
=\infty\})=0$ and $\nu(\{x\in L: \#
(\pi^{-1}(x)\cap I^\bn)
=\infty\})=0$.
\end{enumerate}
Then, by Theorem 1.4.5 in \cite{kig}, 
$\mu(K_w)=N^{-|w|}$ for every $w\in \bigcup_{n\in\bz_+}S^n$  and $\nu(L_w)=N_I^{-|w|}$ for every $w\in \bigcup_{n\in\bz_+}I^n$.
It also holds that $\mu(L)=0$.

Suppose that we are given a strong local regular 
Dirichlet form $(\ce,\cf)$ on $L^2(K,\mu)$. 
$\cf$ is equipped with a norm $\|f\|_{\cf}
=(\ce(f)+\|f\|_{L^2(\mu)}^2)^{1/2}$.
Here and throughout the paper, for each quadratic form $E(\cdot,\cdot)$, 
we abbreviate $E(f,f)$ as $E(f)$.
We assume the following.
\begin{enumerate}
\item[(A5)] (Self-similarity) For each $f\in \cf$ and $i\in S$, $F_i^*f\in \cf$ 
where $F_i^*f=f\circ F_i$. Further, there exists $\rho>0$ such that 
$$\ce (f)=\rho\sum_{i\in S}\ce (F_i^*f),\quad
f\in\cf.$$
\item[(A6)] For every $\Psi\in\fg$, $\Psi^*\cf=\cf$, that is, $\{f\circ \Psi: f\in\cf\}=\cf$.
Further, $\ce(\Psi^* f)=\ce(f)$ for all $f\in\cf$.
\item[(A7)] Let $d_w=(\log \rho N)/(\log \alpha)$. Then $d_w>d_f-d$.
\item[(B1)] 
The space $\cf$ is compactly imbedded in $L^2(K,\mu)$, and $\ce(f)=0$ if and only if $f$ is a constant function.
\end{enumerate}
For each subset $A$ of $W^m$ for some $m\in\bz_+$,
let $\cf_A$ be a function space on $K_A$ such that $\{f|_{K_A}:f\in\cf\}\subset \cf_A\subset \{f\in L^2(K_A): 
F_w^* f\in\cf \mbox{ for all }w\in A\}$.
The space $\cf_A$ will be specified later 
for some class of Dirichlet forms in Section 4. 
Define, for $f,g\in\cf_A$, 
\begin{equation}\label{eq:energy}
\ce_A (f,g)=\rho^m\sum_{w\in A}\ce (F_w^*f,F_w^*g).
\end{equation}
We assume that $\cf_A=\cf_{A\cdot S^n}$ for all $n\in\bn$ and $(\ce_A,\cf_A)$ is a closed form on 
$L^2(K_A,\mu|_{K_A})$.
In what follows, we always consider $\cf_A$ as a normed space with norm $\|f\|_{\cf_A}=(\ce_A(f)+\|f\|_{L^2(K_A)}^2)^{1/2}$.
Due to (A5), $\ce_A(f)=\ce_{A\cdot S^n}(f)$ holds for any $f\in\cf_A$, and $\ce_{\Phi(A)}(f)=\ce_A(f)$ if $\#\Phi(A)=\# A$ by (A6).
When $A=\{w\}$, we use the notation $\ce_w$ in place of $\ce_{\{w\}}$.
Functions in $\cf$ can be naturally considered as elements in $\cf_A$ by the restriction of the domain.
We often write simply $f$ in place of $f|_{K_A}$ when we regard $f\in\cf$ as an element of $\cf_A$, for notational conveniences.
\begin{defn}\label{defn:good}
Let $A$ be a nonempty subset of $W^m$ for some $m\in\bz_+$. 
We say that $A$ is $\ce_A$-connected if, for $f\in\cf_A$, $\ce_A(f)=0$ implies that $f$ is constant on $K_A$.
\end{defn}
\begin{defn}\label{def:same}
Let $A\subset W^m$ and $B\subset W^n$ for some $m$ and $n$.
We say that $A$ and $B$ are of the same type if there exist a homeomorphism $F:K_A\to K_B$ and a bijection $\chi:A\to B$ such that 
$F\circ F_u=F_{\chi(u)}$ for all $u\in A$ and $F^*(\cf_B)=\cf_A$.
\end{defn}

We assume the following.
\begin{enumerate}
\item[(B2)] There exists $\hat I\subset W$ such that the following hold.
\begin{enumerate}
\item[(1)] $\hat I\supset I$ and $\# \hat I<N$.
\item[(2)] For each $w\in I^n$, $\cn_M(w)\cap \hat I^n$ is an $\ce_{\cn_M(w)\cap \hat I^n}$-connected set.
\item[(3)] 
There exist finite elements $u_1,\ldots,u_k\in\bigcup_{n\in\bz_+}I^n$ such that, for any $w\in\bigcup_{n\in\bz_+}I^n$, there exists $j\in\{1,\ldots,k\}$ such that $\cn_M(w)\cap\hat I^{|w|}$ and $\cn_M(u_j)\cap\hat I^{|u_j|}$ are of the same type, and moreover, $F(L_{\cn_M(w)\cap\hat I^{|w|}})=L_{\cn_M(u_j)\cap\hat I^{|u_j|}}$ where $F$ is provided in Definition~\ref{def:same}.
\item[(4)] $C_1:=\sup_{n\in\bz_+}\max_{w\in \hat I^n}\#(\cn_M(w)\cap I^n)<\infty$ and $C_2:=\sup_{n\in\bz_+}\max_{w\in S^n}\#\{v\in\hat I^n: \Phi(v)=w\}<\infty$.
\end{enumerate}
\end{enumerate}
For an open set $U\subset K$, define the capacity of $U$ by
\[\Cp(U)=\inf \{ \|u\|_\cf^2: u\in \cf, u\ge 1
\mbox { $\mu$-a.e.\ on } U\}.\]
The capacity of any set $D\subset K$ is defined as the infimum of the capacity of
open sets that contain $D$.
We denote a quasi-continuous modification of $f\in\cf$ by $\tilde f$.
We assume the following.
\begin{enumerate}
\item[(A8)] There exists some $c>0$ such that $\nu(D)\le c\Cp(D)$ for every compact set $D\subset K$.
\end{enumerate}
By Theorem~3.1 of \cite{ben}, (A8) is equivalent to the following.
\begin{enumerate}
\item[(A8)'] The measure $\nu$ 
charges no set of zero capacity 
and $f\mapsto \tilde f|_L$ is a continuous map from $\cf$ to  $L^2(L,\nu)$.
\end{enumerate}
We will provide sufficient conditions for (A8) in Section~4.

For each $n\in\bz_+$, define $Q_n: L^1(L,\nu)\to \br^{I^n}$ as 
\[Q_nf(w)=\mint_{L_w}f(y)d\nu (y),~~w\in I^n,\]
where in general $-\!\!\!\!\!\int_A \cdots d\lambda(y):=\lambda(A)^{-1}\int_A \cdots d\lambda(y)$ denotes
the normalized integral on $A$. 
Then, one can easily check 
\begin{equation}\label{eq:qmqm1}
N_I^{-1}\sum_{j\in I} Q_{m+1}f(w\cdot j)=Q_mf(w),~~w\in I^m.\end{equation}
Let $m\in\bn$, $A\subset S^m$, and $J\subset I^m$. Define
$\cf(J,A)=\{f\in \cf: f=0 \mbox{ on }K_{S^m\setminus A},\ Q_m (\tilde f|_{L})=0 
\mbox{ on }J\}$,
and define a closed subspace $\ch(J,A)$ of $\cf$ by
\begin{eqnarray*}
\ch(J,A)=\{h\in \cf:  \ce(h,f)=0~~\mbox {for all }
f\in \cf(J,A)\}.
\end{eqnarray*}
When $J$ is an empty set, we omit it from the notation.
We assume the following.
\begin{enumerate}
\item[(B3)] There exist some $l_0,m_0\in\bz_+$, $C>0$, a proper subset $D'(w)$ of $S^{|w|}$ with $w\in D'(w)$
for each $w\in\bigcup_{n\in\bz_+}\Phi(\hat I^{n+m_0})$,  
a finite subset $\Xi\subset \bigcup_{n\in\bz_+}\Phi(\hat I^{n+m_0})$ 
and subsets $D^\sharp(v)$ of $D'(v)$ with $v\in D^\sharp(v)$ for each $v\in \Xi$
such that the following hold.
\begin{enumerate}
\item[(1)] For each $w\in\bigcup_{n\in\bz_+}\Phi(\hat I^{n+m_0})$,
\begin{enumerate}
\item[(a)]
  $w\in D'(w)$ and $D'(w)\subset{\cn_{l_0}(w)}\cap (\Phi(\hat I^n)\cdot S^{m_0})$,
\item[(b)]
  there exists $v\in\Xi$ such that 
\begin{eqnarray*}
\lefteqn{F_w^*(\{h\in\ch(I^{|w|},D'(w)):\int_{K_{D'(w)}}h\,d\mu=0,\ \rho^{-|w|}\ce_{D'(w)}(h)\le 1\})}&&\\
&&\subset F_{v}^*(\{h\in\ch(I^{|v|},D^\sharp(v)):\|h\|_{\cf_{D^\sharp(v)}}\le C\}).
\phantom{\hspace{10em}}
\end{eqnarray*}
\end{enumerate}
\item[(2)] For each $v\in\Xi$,
the operator $F_v^*:\ch(D^\sharp(v))|_{K_{D^\sharp(v)}}\to \cf$ is a compact operator, where $\ch(D^\sharp(v))|_{K_{D^\sharp(v)}}$ is regarded as a subspace of $\cf_{D^\sharp(v)}$.
\end{enumerate}
\end{enumerate}
We set $D(w)=D'(\Phi(w))$ for $w\in\bigcup_{n\in\bz_+}\hat I^{n+m_0}$.
We have a sufficient condition concerning (B3); see Section~4.

The following assumption (B4) will be used in the restriction theorem.
\begin{enumerate}
\item[(B4)]
For $f\in\cf$, if $\ce_{S^m\setminus \Phi(\hat I^m)}(f)=0$ for every $m\in\bz_+$, then $f$ is a constant function.
\end{enumerate}

We next introduce Besov spaces.
\begin{defn}\label{defn:Besov}
For $1\le p<\infty,~1\le q\le \infty,~\beta\ge 0$ and
$m\in \bz_+$, set
\[
a_m (\beta, f):=\gamma^{m\beta}\left(\gamma^{md_f}
\int\!\!\int_{\{(x,y)\in K\times K:\sd(x,y)<c\gamma^{-m}\}}|f(x)-f(y)|^p\,d\mu (x)
d \mu(y)\right)^{1/p}
\]
for $f\in L^p(K,\mu)$,
where $1<\gamma<\infty,~0<c<\infty$. Define a {\it Besov space}
$\Lambda^\beta_{p,q}(K)$ as a set of all $f\in L^p(K,\mu)$ such that
${\bar a}(\beta, f):=\{a_m(\beta, f)\}_{m=0}^{\infty}\in l^q$. 
$\Lambda^\beta_{p,q}(K)$ is a Banach space with the 
norm $\lll f\lll_{\Lambda^\beta_{p,q}(K)}
:=\|f\|_{L^p(K)}+\|{\bar a}(\beta, f)\|_{l^q}$. 
Let $\hat\Lambda^\beta_{p,q}(K)$ denote the closure of $\Lambda^\beta_{p,q}(K)\cap C(K)$ in $\Lambda^\beta_{p,q}(K)$.
$\Lambda^\beta_{p,q}(L)$ and $\hat\Lambda^\beta_{p,q}(L)$ are defined in the same way by replacing $(K,\mu)$ by $(L,\nu)$.
\end{defn}
We remark that this definition is valid for general Alfors regular compact sets $K$ with normalized Hausdorff measure $\mu$.
We use the notation $\Lambda^\beta_{p,q}(K)$ following \cite{gri}.
$\Lambda^\beta_{p,q}(K)$ was denoted by Lip $(\beta,p,q)(K)$ in \cite{j2,kum} 
and by $\Lambda_\beta^{p,q}(K)$ in \cite{stri}. 
Note that different choices of $c>0$ and $\gamma>1$ provide the same space $\Lambda^\beta_{p,q}(K)$ with equivalent norms. In what follows, we will take $\gamma=\alpha$.

We are now ready to state our main theorems. 
Let $\beta=d_w/2-(d_f-d)/2$.
\begin{theo}\label{theo:mainthm1}
   Suppose that \rom{(A1)}--\rom{(A8)} and \rom{(B1)}--\rom{(B4)} hold.
   Then, for every $f\in\cf$, $\tilde f|_L$ belongs to $\hat\Lambda^\beta_{2,2}(L)$.
   Moreover, there exists $c>0$ such that $\|\tilde f|_L\|_{\Lambda^\beta_{2,2}(L)}\le c\|f\|_{\cf}$ for every $f\in\cf$.
\end{theo}
\begin{theo}\label{theo:mainthm2}
  Suppose that \rom{(A1)}--\rom{(A8)} and \rom{(C1)}--\rom{(C2)} hold. (The conditions \rom{(C1)} and \rom{(C2)} will be defined in Section~3.3). 
  Then, there exists a bounded linear map $\xi$ from $\hat\Lambda^\beta_{2,2}(L)$ to $\cf$ such that $\xi(\Lambda^\beta_{2,2}(L)\cap C(L))\subset\cf\cap C(K)$ and $\widetilde{\xi f}|_L=f$ $\nu$-a.e.\ for all $f\in\hat\Lambda^\beta_{2,2}(L)$.
\end{theo}
In what follows, we often write $\cf|_L=\hat\Lambda^\beta_{2,2}(L)$ to denote the assertions of two theorems above.
\begin{rem}\label{rem:2.5} 
In the following two cases, we can prove 
$\hat\Lambda^\beta_{2,2}(L)=\Lambda^\beta_{2,2}(L)$.\\
1) $L\subset \br^n$ for some $n\in \bn$ and $\beta<1$.~~In this case, the following trace theorem 
holds due to \cite{jw}; $B^{2,2}_{\beta+(n-d)/2}(\br^n)|_L=\Lambda^\beta_{2,2}(L)$ 
where $B^{2,2}_{\gamma}(\br^n)$ is the classical Besov space with smoothness order $\gamma$. Since 
$C_0^\infty(\br^n)$ is dense 
in $B^{2,2}_{\gamma}(\br^n)$ for $\gamma>0$, it follows that functions from 
$C_0^\infty(\br^n)$ restricted to $L$ are dense in $\Lambda^\beta_{2,2}(L)$.\\
2) $\beta>d/2$.~~
In this case, the following holds due to \cite{gri} Theorem 8.1;
$\Lambda_{2,\infty}^\beta(L)\subset {\cal C}^{\beta-d/2}(L)$, where
${\cal C}^{\lambda}(L)$ is a H\"older space defined as follows. $u\in {\cal C}^{\lambda}(L)$ if 
\begin{equation}
\|u\|_{{\cal C}^\lambda(L)}:=\|u\|_{L^\infty(L)}+
\mathop{\nu\rm -esssup}_{x,y\in L,\ x\ne y}
\frac{|u(x)-u(y)|}{\sd(x,y)^\lambda}<\infty.
\end{equation}
Since $\Lambda_{2,2}^\beta(L)\subset \Lambda_{2,\infty}^\beta(L)$, we see that any element in 
$\Lambda_{2,2}^\beta(L)$ is continuous in this case. 
\end{rem}
\begin{rem}\label{rem:timechange}
Since $\nu$ is smooth with respect to $(\ce,\cf)$, we can consider the time changed Markov process with respect to the positive continuous additive functional associated with $\nu$ via the Revuz correspondence.
By the general theory of Dirichlet forms, this has an associated regular Dirichlet form $(\check\ce,\check\cf)$ on $L^2(L,\nu)$ with
$\check\cf=\{f\in L^2(L,\nu): f=\tilde u \ \nu\mbox{-a.e.\ on $L$ for some }u\in\cf_e\}$,
where $\cf_e$ is the family of $\mu$-measurable functions $u$ on $K$ such that $|u|<\infty$ $\mu$-a.e.\ and there exists an $\ce$-Cauchy sequence $\{u_n\}_{n\in\bn}$ of functions in $\cf$ such that $\lim_{n\to\infty}u_n=u$ $\mu$-a.e.
As is seen in the proposition below, $\cf_e=\cf$ in our framework.
So, our main theorems determine the function space $\check\cf$.
\end{rem}
\begin{pro}
  Under the condition (B1), $\cf_e=\cf$.
\end{pro}
\Proof
  By (B1), there exists some $c>0$ such that 
  \begin{equation}\label{eq:poincare}
  \left\|f-\int_K f\,d\mu\right\|_{L^2(K)}^2\le c\ce(f),
  \quad f\in\cf.
  \end{equation}
  Let $u\in\cf_e$. Take $\{u_n\}_{n\in\bn}$ from $\cf$ as in the definition of $\cf_e$ in Remark~\ref{rem:timechange}.
  Define $g_n=u_n-\int_K u_n\,d\mu$ for each $n$.
  Then, $\{g_n\}_{n\in\bn}$ is $\ce$-Cauchy.
  Since $\int_K g_n\,d\mu=0$, (\ref{eq:poincare}) implies that $\{g_n\}$ is also $L^2(K)$-Cauchy.
  Therefore, $g_n$ converges to some $g$ in $\cf$.
  By taking a subsequence, we may assume that $g_n\to g$ $\mu$-a.e.
  Thus, $\int_K u_n\,d\mu\,({}=u_n-g_n)$ converges to some $C\in\br$.
  In particular, $\int_K u_n\,d\mu$ converges to $C$ in $\cf$ as a sequence of constant functions.
  Therefore, $u_n$ converges to $g+C$ in $\cf$.
  This implies that $u=g+C$ belongs to $\cf$.
\QED
\section{Proof of main theorems}
\subsection{Discrete approximation} 
In this section, we assume (A1)--(A8).
For $n\in\bz_+$, define a bilinear form on $I^n$ as 
\[E_{(n)}(g,g)=\sum_{v,w\in I^n,~ v\stackrel{n,L}{\lr} w}
(g(v)-g(w))^2~~~\mbox{for }~ g\in \br^{I^n}.\]
We then have the following discrete characterization of $ \Lambda^\beta_{2,q}(L)$ (for related results,
see \cite{kam}). 

\begin{lem}\label{theo:lem2}
Let  $\beta>0$ and $q\in[1,\infty]$. 
Then, there exists $c_1>0$ such that 
for each $f\in L^2(L,\nu)$, 
\begin{eqnarray}
\lefteqn{c_1\left\|\left\{
\alpha^{n\beta}\left(\alpha^{nd}\int\!\!\!\int_{\{(x,y)\in L\times L: 
\sd(x,y)< k_1 \alpha^{-n}\}}|f(x)-f(y)|^2\,d\nu(x)\,d\nu(y)
\right)^{1/2}\right\}_{n=0}^\infty\right\|_{l^q}}\nonumber\\
&\le&\left\|\left\{\alpha^{n\beta}\left(\alpha^{-nd}E_{(n)}(Q_n f)
\right)^{1/2}\right\}_{n=0}^{\infty}\right\|_{l^q}\nonumber\\
&\le& \left\|\left\{
\alpha^{n\beta}\left(\alpha^{nd}\int\!\!\!\int_{\{(x,y)\in L\times L: 
\sd(x,y)< k_2 \alpha^{-n}\}}|f(x)-f(y)|^2\,d\nu(x)\,d\nu(y)\right)^{1/2}
\right\}_{n=0}^\infty\right\|_{l^q}.\label{eq:equivalence}
\end{eqnarray}
Here, $k_1$ and $k_2$ are provided in \rom{(A2)}.
\end{lem}
\Proof
Due to the choice of $M$, the exists some $c_2>0$ such that
\[
\sum_{i\in I}\left(g(i)-N_I^{-1}
\sum_{j\in I}g(j)\right)^2\le c_2 E_{(1)}(g),\quad g\in \br^I.
\]
For $f\in L^2(L,\nu)$ and $n\in\bz_+$, we have
\begin{eqnarray*}
\lefteqn{\int\!\!\!\int_{\{(x,y)\in L\times L: \sd(x,y)< k_1 
\alpha^{-n}\}}|f(x)-f(y)|^2\,d\nu(x)\,d\nu(y)}\\
&\le& \sum_{(v,w)\in I^n\times I^n,\,v\stackrel{n,L}{\lr} w}\int\!\!\!
\int_{L_v\times L_w}|f(x)-f(y)|^2\,d\nu(x)\,d\nu(y)
\qquad\mbox{(by (A2))}\\
&\le& \sum_{(v,w)\in I^n\times I^n,\,v\stackrel{n,L}{\lr} w}\int\!\!\!
\int_{L_v\times L_w}3\{|f(x)-Q_nf(v)|^2+|Q_nf(v)-Q_nf(w)|^2\\
&&{}+|Q_nf(w)-f(y)|^2\}\,d\nu(x)\,d\nu(y)\\
&\le& 6C_0 N_I^{-n}\sum_{v\in I^n} \int_{L_v} (f(x)-Q_nf(v))^2\,d\nu(x)
+3 N_I^{-2n}E_{(n)}(Q_n f),
\end{eqnarray*}
where $C_0$ is what appeared in (A1).
Concerning the first term, we have
\begin{eqnarray*}
\lefteqn{\sum_{v\in I^n} \int_{L_v} (f(x)-Q_nf(v))^2\,d\nu(x)}\\
&=& \int_L f(x)^2\,d\nu(x)-N_I^{-n}\sum_{v\in I^n} Q_nf(v)^2\\
&=&\sum_{m=n}^\infty \left(N_I^{-(m+1)}\sum_{v\in I^{m+1}} 
Q_{m+1}f(v)^2-N_I^{-m}\sum_{w\in I^m} Q_mf(w)^2\right)\\
&=&\sum_{m=n}^\infty N_I^{-(m+1)}\sum_{w\in I^m}\sum_{i\in I}
\left(Q_{m+1}f(w\cdot i)-N_I^{-1}\sum_{j\in I}Q_{m+1}f(w\cdot j)\right)^2\\
&\le&c_2\sum_{m=n}^\infty N_I^{-(m+1)}\sum_{w\in I^m}E_{(1)}(Q_{m+1}f(w\cdot*))\\
&\le&c_2\sum_{m=n}^\infty N_I^{-(m+1)}E_{(m+1)}(Q_{m+1}f),
\end{eqnarray*}
where the martingale convergence theorem was used in the second equality and (\ref{eq:qmqm1}) was used in the third equality.
Note that $\alpha^d=N_I$.
Suppose $q\in[1,\infty)$. Then,
\begin{eqnarray*}
\lefteqn{\sum_{n=0}^\infty
\alpha^{n(\beta+d/2)q}\left(\int\!\!\!\int_{\{(x,y)\in L\times L: 
\sd(x,y)< k_1 \alpha^{-n}\}}|f(x)-f(y)|^2\,d\nu(x)\,d\nu(y)\right)^{q/2}}\\
&\le&\sum_{n=0}^\infty\alpha^{n(\beta+d/2)q}\left(6c_2C_0
N_I^{-n}\sum_{m=n}^\infty N_I^{-(m+1)}E_{(m+1)}(Q_{m+1}f)
+3 N_I^{-2n}E_{(n)}(Q_n f)\right)^{q/2}\\
&\le&c_3\sum_{n=0}^\infty\alpha^{n\beta q}\left(\sum_{m=n}^\infty
\alpha^{-md}E_{(m)}(Q_{m}f)\right)^{q/2}\\
&\le&c_4\sum_{m=0}^\infty\alpha^{m(\beta-d/2)q}E_{(m)}(Q_{m}f)^{q/2}\\
&=& c_4\left\|\left\{\alpha^{n\beta}\left(\alpha^{-nd}
E_{(n)}(Q_n f)\right)^{1/2}\right\}_{n=0}^\infty\right\|_{l^q}^q,
\end{eqnarray*}
where in the third inequality, we used (A7) and the following inequality for $a>0$:
\begin{equation}\label{eq:harheld}
\sum_{i=0}^\infty 2^{ai}\left(\sum_{j\in \Lambda_i} a_j\right)^p
\le c \sum_{j=0}^\infty 2^{aj}a_j^p
\quad\mbox{for }a\ne 0,\ p>0,\ a_j\ge0,
\end{equation}
where $\Lambda_i=\{i,i+1,\ldots\}$ when $a>0$ and $\Lambda_i=\{0,1,\ldots, i\}$ when $a<0$.
When $0<p\le1$, this is obvious since $(x+y)^p\le x^p+y^p$ for $x$, $y\ge0$.
When $p>1$, this is proved by applications of H\"older's inequality; see e.g.\ \cite{lei}.

When $q=\infty$, letting $\gamma=\left\|\left\{\alpha^{n\beta}\left(\alpha^{-nd}
E_{(n)}(Q_n f)\right)^{1/2}\right\}_{n=0}^\infty\right\|_{l^\infty}$, we have for every
$n\in\bz_+$,
\begin{eqnarray*}
\lefteqn{\left|
\alpha^{n(\beta+d/2)}\left(\int\!\!\!\int_{\{(x,y)\in L\times L: \sd(x,y)< k_1 
\alpha^{-n}\}}|f(x)-f(y)|^2\,d\nu(x)\,d\nu(y)\right)^{1/2}\right|^2}\\
&\le&c_5\alpha^{n(2\beta+d)}\left(N_I^{-n}\sum_{m=n}^\infty N_I^{-(m+1)}
E_{(m+1)}(Q_{m+1}f)+N_I^{-2n}E_{(n)}(Q_n f)\right)\\
&\le&c_5\alpha^{2n\beta}\sum_{m=n}^\infty \alpha^{-2m\beta}\gamma^2\\
&=&\frac{c_5}{1-\alpha^{-2\beta}}\gamma^2.
\end{eqnarray*}
Thus, the first inequality in (\ref{eq:equivalence}) is proved.

Next, we have
\begin{eqnarray*}
E_{(n)}(Q_n f)
&=&\sum_{(v,w)\in I^n\times I^n,\,v\stackrel{n,L}{\lr} w} \left|N_I^{2n}\int\!\!\!
\int_{L_v\times L_w}\{f(x)-f(y)\}\,d\nu(x)\,d\nu(y)\right|^2\\
&\le&\sum_{(v,w)\in I^n\times I^n,\,v\stackrel{n,L}{\lr} w} N_I^{2n}\int\!\!\!
\int_{L_v\times L_w}|f(x)-f(y)|^2\,d\nu(x)\,d\nu(y)\\
&\le&\alpha^{2nd}\int\!\!\!\int_{\{(x,y)\in L\times L: \sd(x,y)< 
k_2 \alpha^{-n}\}}|f(x)-f(y)|^2\,d\nu(x)\,d\nu(y),
\end{eqnarray*}
which deduces the second inequality of (\ref{eq:equivalence}).
\QED
\begin{rem} Quite recently, M. Bodin (\cite{bod}) gives a discrete characterization of 
$\Lambda_{p,q}^\beta(K)$ for the Alfors $d$-regular set $K$ if it has a 
regular triangular system with some property (property~(B) in the thesis). 

\end{rem}
\subsection{Proof of the restriction theorem}
In this section, we assume (A1)--(A8) and (B1)--(B4), and prove Theorem~\ref{theo:mainthm1}. 
The following lemma is immediately proved by equation~(\ref{eq:energy}).
\begin{lem}\label{lem:scale}
  Let $A\subset W^m$, $B\subset W^n$, $f\in\cf_A$, and $g\in\cf_B$.
  Suppose that there exists a bijection $\iota$ from $A$ to $B$ and $F_v^* f= F_{\iota(v)}^* g$ for every $v\in A$. Then, $\rho^{-m}\ce_A(f)=\rho^{-n}\ce_B(g)$.
\end{lem}
  Let $n\in\bz_+$ and $w\in I^n$. 
  Let $A=\cn_M(w)\cap \hat I^n$.
  Define $\cg_w=\{f\in\cf_A: Q_n (\tilde f|_{L_A})=0 \mbox{ on }\cn_M(w)\cap I^n\}$ and $\ck_w=
  \{h\in\cf_A: \ce_A(h,f)=0 \mbox{ for all }f\in\cg_w\}$.
  Here, we used (and will use) notations $Q_n (\tilde f|_{L_A})$ (on $A$) and $\ce_A(f)$ 
  for $f\in \cf_A$ in the obvious sense.
\begin{lem}\label{lem:pre}
\begin{enumerate}
\item[$(1)$] There exists some $c>0$ such that $\|f\|_{L^2(K_A)}^2\le c\ce_A(f)$ for all $f\in\cg_w$.
\item[$(2)$] For each $g\in\cf_A$, there exists $h_g\in \ck_w$ such that $Q_n (\tilde h_g|_{L_A})=Q_n (\tilde g|_{L_A})$ on $\cn_M(w)\cap I^n$ and $\ce_A(h_g)\le \ce_A(g)$.
\end{enumerate}
\end{lem}
\Proof
  (1) Suppose that the claim does not hold.
  Then, there exists a sequence $\{f_k\}_{k\in\bn}\subset\cg_w$ such that $\|f_k\|_{L^2(K_A)}=1$ and $\lim_{k\to\infty}\ce_A(f_k)=0$.
  We may assume that $f_k$ converges weakly to some $f$ in $\cf_A$ and $F_w^* f_k$ converges to $F_w^* f$ weakly in $\cf$ for every $w\in A$.
  By (B1), $F_w^* f_k$ converges to $F_w^* f$ in $L^2(K)$ for each $w\in A$. Thus, $f_k$ converges to $f$ in $L^2(K_A)$.
  We also have
  $\ce_A(f)\le \liminf_{k\to\infty}\ce_A(f_k)=0$.
  Therefore, $\ce_A(f)=0$.
  In view of (B2)(2), $f$ is constant on $K_A$.
  Since $f$ belongs to $\cg_w$ by (A8)', we conclude $f=0$ on $K_A$, which is a contradiction to the fact that $\|f\|_{L^2(K_A)}=\lim_{k\to\infty}\|f_k\|_{L^2(K_A)}=1$.
  
  (2) Let $\cf_g=\{f\in\cf_A:Q_n (\tilde f|_{L_A})=Q_n(\tilde g|_{L_A})\mbox{ on }\cn_M(w)\cap I^n\}$.
  Take a sequence $\{h_k\}_{k\in\bn}\subset \cf_g$ such that $\ce_A(h_k)$ converges to the infimum of $\{\ce_A(f):f\in\cf_g\}$.
  Since
  \begin{equation}\label{eq:domination}
    \|h_k\|_{L^2(K_A)}
    \le \|h_k-g\|_{L^2(K_A)}
       +\|g\|_{L^2(K_A)}
    \le c^{1/2}\ce_A(h_k-g)^{1/2}
       +\|g\|_{L^2(K_A)},
  \end{equation}
  we have $\sup_k \|h_k\|_{L^2(K_A)}<\infty$.
  There exists a weak limit $h\in\cf_A$ of a subsequence of $\{h_k\}_{k\in\bn}$ in $\cf_A$.
  Then $h\in\cf_g$ and $h$ attains the infimum of $\{\ce_A(f): f\in\cf_g\}$.
  Dividing by $\epsilon$ both sides of the inequality $\ce_A(h+\epsilon f)-\ce_A(h)\ge 0$ for $f\in\cg_w$ and letting $\epsilon\to0$, we obtain $h\in\ck_w$.
\QED
\begin{lem}\label{lem:qnapp}
There exists some $c_1>0$ such that
\begin{equation}\label{eq:qnapp}
c_1\rho^{-n}\ce_{\Phi(\hat I^n)} (f)\ge E_{(n)}(Q_n(\tilde f|_L))~~~\mbox{for all $f\in \cf$ and $n\in\bz_+$}.
\end{equation}
\end{lem}
\Proof
  First, we prove that $\ck_w$ is a finite dimensional vector space.
  For each $i\in\cn_M(w)\cap I^n$, take a function $g_i\in\cf$ such that
  $Q_n(\tilde g_i|_L)(j)=\left\{\begin{array}{ll}1& j=i\\0& j\ne i\end{array}\right.$ for all $j\in\cn_M(w)\cap I^n$.
  Existence of such functions is assured by the regularity of the Dirichlet form $(\ce,\cf)$.
  Define a linear map $\Theta\colon \ck_w\to\br^{\cn_M(w)\cap I^n}$ by $\Theta(f)=\{\ce_A(f,g_i)\}_{i\in\cn_M(w)\cap I^n}$.
  Suppose $f$ belongs to the kernel of $\Theta$.
  Then $\ce_A(f,g)=0$ for every $g\in\cf_A$, which implies that $f$ is constant on $K_A$ by (B2)~(2).
  Therefore, $\ck_w$ is finite dimensional.
  
Since $\ce_{A}(h)=0$ implies
$\sum_{v\in A}(Q_n (\tilde h|_{L_A})(v)-Q_n (\tilde h|_{L_A})(w))^2=0$ for $h\in\ck_w$, there exists $c_2>0$ such
that $ \sum_{v\in A}(Q_n (\tilde h|_{L_A})(v)-Q_n (\tilde h|_{L_A})(w))^2\le c_2\rho^{-n}\ce_A (h)$ for every $h\in\ck_w$.
  By (B2)~(3) and Lemma~\ref{lem:scale}, we can take $c_2$ independently with respect to $w\in \bigcup_{n\in\bz_+}I^n$.
  Therefore, for any $f\in \cf$ and $n\in\bz_+$,
  by taking $h_f\in\ck_w$ as in Lemma~\ref{lem:pre} (2),
  \begin{eqnarray*}
  \sum_{v\in \cn_M(w)\cap I^n}(Q_n (\tilde f|_L)(v)-Q_n(\tilde f|_L)(w))^2
  &=&\sum_{v\in \cn_M(w)\cap I^n}(Q_n(\tilde h_f|_{L_A})(v)-Q_n(\tilde h_f|_{L_A})(w))^2\\
  &\le&c_2\rho^{-n}\ce_{A} (h_f)\\
  &\le&c_2\rho^{-n}\ce_{A} (f).
  \end{eqnarray*}
  This implies that
  \begin{eqnarray*}
  E_{(n)}(Q_n(\tilde f|_L))
  &=& \sum_{w\in I^n}\sum_{v\in \cn_M(w)\cap I^n}(Q_n (\tilde f|_L)(v)-Q_n(\tilde f|_L)(w))^2\\
  &\le& c_2\rho^{-n}\sum_{w\in I^n}\ce_{\cn_M(w)\cap \hat I^n} (f)\\
  &\le& c_2C_1\rho^{-n}\ce_{\hat I^n} (f)
  \le c_2C_1C_2\rho^{-n}\ce_{\Phi(\hat I^n)} (f),
  \end{eqnarray*}
  where $C_1$ and $C_2$ are provided in (B2) (4).
\QED
Recall finite sets $\Xi$ and $D^\sharp(v)$ for $v\in\Xi$ introduced in (B3).
\begin{lem}\label{lem:compact}
For each $v\in\Xi$,
the operator
$F_v^*\colon \ch(I^{|v|},D^\sharp(v))|_{K_{D^\sharp(v)}}\to \cf$ is a compact operator.
Here, $\ch(I^{|v|},D^\sharp(v))|_{K_{D^\sharp(v)}}$ is regarded as a subspace of $\cf_{K_{D^\sharp(v)}}$.
\end{lem}
\Proof
Define $I(v)=\{w\in I^{|v|}: L_w\not\subset K_{S^{|v|}\setminus D^\sharp(v)}\}$.
Note that $\ch(I^{|v|},D^\sharp(v))=\ch(I(v),D^\sharp(v))$.
For each $i\in I(v)$, take a function $g_i$ in $\cf(D^\sharp(v))$ such that $Q_{|v|}(\tilde g_i|_L)(j)=\left\{\begin{array}{cl}1&\mbox{if }j=i\\0&\mbox{if }j\ne i\end{array}\right.$ for all $j\in I(v)$. 
Define a linear map $\Theta: \ch(I^{|v|},D^\sharp(v))|_{K_{D^\sharp(v)}}\to \br^{I(v)}$ by
$\Theta(f)=\{\ce(f,g_i)\}_{i\in I(v)}$.
Then, the kernel of $\Theta$ is equal to $\ch(D^\sharp(v))|_{K_{D^\sharp(v)}}$.
The homomorphism theorem implies that $\left.\ch(I^{|v|},D^\sharp(v))|_{K_{D^\sharp(v)}}\right/\ch(D^\sharp(v))|_{K_{D^\sharp(v)}}\simeq \Theta(\ch(I^{|v|},D^\sharp(v))|_{K_{D^\sharp(v)}})$ as a vector space.
Therefore, there exists a finite dimensional vector space $Z$ of $\ch(I^{|v|},D^\sharp(v))|_{K_{D^\sharp(v)}}$ such that $\ch(I^{|v|},D^\sharp(v))|_{K_{D^\sharp(v)}}$ is a direct sum of $\ch(D^\sharp(v))|_{K_{D^\sharp(v)}}$ and $Z$.
Condition (B3)~(2) concludes the assertion.
\QED
\begin{lem}\label{lem:harmext}
Let $m\in\bn$, $A$ a proper subset of $S^m$, and $J$ a subset of $I^m$.
For $g\in\cf$, there exists a unique function $g'$ in $\ch(J,A)$ such that $g'=g$ on $K_{S^m\setminus A}$ and $Q_m(\widetilde{g'}|_L)=Q_m(\tilde g|_L)$ on $J$.
  Moreover, 
  there exists $c>0$ such that 
\begin{equation}\label{eq:dominationA}
\|g'\|_{\cf_A}\le c\|g\|_{\cf_A}, \quad
\ce(g')\le \ce(g)
\end{equation}
for all $g\in\cf$.
Further, if $g\ge0$ $\mu$-a.e., then $g'\ge0$ $\mu$-a.e.
\end{lem}
\Proof
  First, we prove that there exists some $c'>0$ such that $\|f\|_{L^2(K_A)}^2\le c'\ce_A(f)$ for every $f\in \cf(A)$.
  Suppose this does not hold. Then, there exists a sequence $\{f_n\}_{n\in\bn}\subset \cf(A)$ such that $\|f_n\|_{L^2(K_A)}=1$ for every $n$ and $\ce_A(f_n)$ converges to $0$ as $n\to\infty$. 
We may assume that
$f_n$ converges weakly to some $f$ in $\cf$. 
Then, $f_n$ converges to $f\in\cf$ in $L^2(K)$ by (B1), and $\ce(f)\le \liminf_{n\to\infty}\ce(f_n)=0$.
Therefore, $\ce(f)=0$ and $f$ is constant on $K$.
Since $f\in\cf(A)$ and $A\ne S^m$, $f$ is identically $0$, which is 
contradictory to the fact $\|f\|_{L^2(K)}=1$.
  
  Now, given $g\in\cf$, let $\cf_g=\{f\in\cf: f=g\mbox{ on }K_{S^m\setminus A} \mbox{ and }Q_m(\tilde f|_L)=Q_m(\tilde g|_L) \mbox{ on }J\}$.
Then, in exactly the same way as the proof of Lemma~\ref{lem:pre}~(2), there exists $h\in\cf_g$ attaining the
infimum of $\{\ce(f):f\in \cf_g\}$ and $h\in \ch(J,A)$. 
Such functions exist uniquely; indeed, if both $h$ and $h'$
attain the infimum above, we have
\[
  \ce\left(\frac{h-h'}2\right)=\frac12\left(\ce(h)+\ce(h')\right)-\ce\left(\frac{h+h'}2\right)\le0,
\]
which implies that $h-h'$ is a
constant. Since $h-h'=0$ on $K_{S^m\setminus A}$, we conclude that $h=h'$.
On the other hand, it is easy to see that $g'$ should attain the infimum above. Therefore, $g'$ is uniquely determined. 
By the inequality similar to (\ref{eq:domination}),
we conclude (\ref{eq:dominationA}).
The last assertion follows from the characterization of $g'$ above and the Markov property of the Dirichlet form.
\QED

The following is the key proposition. Condition (B4) will be used (only) here. 
\begin{pro}\label{theo:keyhin}
There exist $0<c_0<1$ and $b_0\in\bn$ such that the following holds for all 
$n\in\bz_+$ and $h\in \ch(I^{n},\Phi(\hat I^n))$: 
\[
\ce_{\Phi(\hat I^{n+b_0})}(h)\le c_0\ce_{\Phi(\hat I^{n})}(h).
\]
Moreover,
for all $i\ge j\ge 1$, $b=0,1,\dots,b_0-1$, and $h\in \ch(I^{b_0j},\Phi(\hat I^{b_0j}))$, 
\[
\ce_{\Phi(\hat I^{b_0i+b})}(h)\le c_0^{i-j}\ce_{\Phi(\hat I^{b_0j+b})}(h).
\]
\end{pro}
\Proof
It is enough to prove the first claim. 
Recall $l_0$ and $m_0$ in condition (B3).
By (B3), $C:=\sup_{n\in\bz_+}\max_{w\in \hat I^{n+m_0}}\# D(w)$ is finite.
Let $n\in\bz_+$ and $w\in \hat I^{n+m_0}$. 
Define 
\begin{eqnarray*}
\cc_w&=&\{F_w^*f: f\in\ch(I^{n+m_0},D(w)),\ \int_{K_{D(w)}}f\,d\mu=0,\ \rho^{-(n+m_0)}\ce_{D(w)}(f)\le1\},\\
\cc&=&\mbox{the closure of }\bigcup_{w\in \bigcup_{n\in\bz_+}\hat I^{n+m_0}}\cc_w \mbox{ in }\cf.
\end{eqnarray*}
Then, $\cc$ is a compact subset in $\cf$ by Lemma~\ref{lem:compact} and (B3). 
Let $\dl=1/(4C^2)$ and define
$\cc(\dl)=\{f\in\cc: \ce(f)\ge\dl\}$.
Since (B4) holds, for each $f\in\cc(\dl)$, there exist $m(f)\in\bn$ 
and $a(f)\in(0,1)$ such that $\ce_{\Phi(\hat I^m)}(f)<a(f)\ce(f)$ for all $m\ge m(f)$. 
By continuity, $\ce_{\Phi(\hat I^m)}(g)<a(f)\ce(g)$ for all $m\ge m(f)$ for any $g$ in some neighborhood of $f$ in
$\cf$. Since $\cc(\dl)$ is compact in $\cf$, there exist $m_1\in\bn$ and $a_1\in(0,1)$ 
such that $\ce_{\Phi(\hat I^{m_1})}(f)<a_1\ce(f)$ for every $f\in\cc(\dl)$.
In particular, 
\begin{equation}\label{eq:dom1}
\ce(f)\le a_2 \ce_{S^{m_1}\setminus \Phi(\hat I^{m_1})}(f),
\quad f\in\cc(\dl)
\end{equation}
with $a_2=(1-a_1)^{-1}>1$.

Now, take $h$ as in the claim of the proposition.
We construct an oriented graph such that the set of vertices is $\Phi(\hat I^{n+m_0})$ and a set of oriented edges is 
$E=\{(v,w)\in \Phi(\hat I^{n+m_0})\times \Phi(\hat I^{n+m_0}): 
v\in D'(w),\, \ce_{w}(h)>0 \mbox{ and }
\ce_{w}(h)\ge 2C \ce_{v}(h)\}$.
This graph does not allow any loops.
Let $Y$ be the set of all elements $w$ in $\Phi(\hat I^{n+m_0})$ such that $\ce_w(h)>0$ and $w$ is not a source of any edges.
For $w\in Y$, define $N_0(w)=\{w\}$, $N_k(w)=\{v\in \Phi(\hat I^{n+m_0})\setminus \bigcup_{l=0}^{k-1}
N_l(w): (v,u)\in E \mbox{ for some }u\in N_{k-1}(w)\}$ for $k\in\bn$ 
inductively, and $N(w)=\bigcup_{k\ge0}N_k(w)$.
It is clear that $\# N_k(w)\le C^k$ and $\ce_{v}(h)\le 
(2C)^{-k}\ce_{w}(h)$ for all $k\ge0$ and $v\in N_k(w)$.
Then, for each $w\in Y$,
\begin{equation}\label{eq:dom2}
  \ce_{N(w)}(h)
  = \sum_{k=0}^\infty\sum_{v\in N_k(w)}\ce_{v}(h)
  \le \sum_{k=0}^\infty C^k (2C)^{-k}\ce_{w}(h)
  = 2\ce_{w}(h).
\end{equation}
Suppose $w\in Y$ and $\ce_w(h)\ge \dl\ce_{D'(w)}(h)$. Then, since 
\[
F_w^*\left(\left(h-\mint_{K_{D'(w)}}h\,d\mu\right)\times \rho^{(n+m_0)/2}\ce_{D'(w)}(h)^{-1/2}\right)\in\cc(\dl),
\]
(\ref{eq:dom1}) implies that $\ce(F_w^* h)\le a_2 \ce_{S^{m_1}\setminus \Phi(\hat I^{m_1})}(F_w^* h)$, namely,
\[
  \ce_{w}(h)\le a_2\ce_{w\cdot (S^{m_1}\setminus \Phi(\hat I^{m_1}))}(h).
\]
Next, suppose $w\in Y$ and $\ce_w(h)< \dl\ce_{D'(w)}(h)$.
Since $w$ is not a source of any edges, $\ce_v(h)<2C\ce_w(h)$ for every $v\in D'(w)\cap \Phi(\hat I^{n+m_0})$.
Then,
\[
  \ce_{D'(w)\cap \Phi(\hat I^{n+m_0})}(h)
  < C\cdot 2C\ce_w(h)
  < 2C^2\dl \ce_{D'(w)}(h)
  = \frac12 \ce_{D'(w)}(h),
\]
which implies 
$\ce_{D'(w)\cap \Phi(\hat I^{n+m_0})}(h)
< \ce_{D'(w)\cap ((\Phi(\hat I^n)\cdot S^{m_0})\setminus \Phi(\hat I^{n+m_0}))}(h)$
by (B3)~(1)~(a).
In particular,
\[
  \ce_w(h)
  <\ce_{D'(w)\cap ((\Phi(\hat I^n)\cdot S^{m_0})\setminus \Phi(\hat I^{n+m_0}))}(h).
\]
Therefore, in any cases, we have for $w\in Y$, 
\begin{eqnarray}\label{eq:dom3}
  \ce_w(h)&\le&a_2\ce_{w\cdot(S^{m_1}\setminus \Phi(\hat I^{m_1}))\cup((D'(w)\cap ((\Phi(\hat I^n)\cdot S^{m_0})\setminus \Phi(\hat I^{n+m_0})))\cdot S^{m_1})}(h) \nonumber\\
  &\le& a_2\ce_{(D'(w)\cdot S^{m_1})\cap((\Phi(\hat I^n)\cdot S^{b_0})\setminus \Phi(\hat I^{n+b_0}))}(h),
\end{eqnarray}
where $b_0=m_0+m_1$ and note that
$\Phi(\hat I^{n+b_0})\subset \Phi(\hat I^{n+m_0})\cdot S^{m_1}$.
Then we have
\begin{eqnarray*}
\ce_{\Phi(\hat I^{n+b_0})}(h)
&\le& \ce_{\Phi(\hat I^{n+m_0})}(h)\\
&\le& \sum_{w\in Y} \ce_{N(w)}(h)\\
&\le& 2a_2\sum_{w\in Y}\ce_{(D'(w)\cdot S^{m_1})\cap((\Phi(\hat I^n)\cdot S^{b_0})\setminus \Phi(\hat I^{n+b_0}))}(h)\\
&\le& 2a_2 C_3\ce_{(\Phi(\hat I^n)\cdot S^{b_0})\setminus \Phi(\hat I^{n+b_0})}(h).
\end{eqnarray*}
Here, we used (\ref{eq:dom2}) and (\ref{eq:dom3}) in the third inequality
and $C_3:=\sup_{n\in\bz_+}\max_{v\in S^{n+m_0}}\#(\cn_{l_0}(v)\cap S^{n+m_0})$ is finite by (A1).
Hence, the claim of the proposition holds with $c_0=2a_2 C_3/(1+2a_2 C_3)$.
\QED
\Proofwithname{Proof of Theorem~\ref{theo:mainthm1}.}
Fix $b\in\{0,1,\ldots,b_0-1\}$.
For each $f\in \cf$, we can take $g_m\in \ch(I^{b_0m+b},\Phi(\hat I^{b_0m+b}))$ such that $g_m=f$ on $K_{S^{b_0m+b}
\setminus \Phi(\hat I^{b_0m+b})}$ and $Q_{b_0m+b}(\tilde g_m|_L)=Q_{b_0m+b}(\tilde f|_L)$ by Lemma~\ref{lem:harmext}.
By using the relations $\|g_m\|_\cf\le c\|f\|_\cf$, $\ce(g_m)\le \ce(f)$ (by Lemma~\ref{lem:harmext}), and $g_m\to f$ $\mu$-a.e., we will prove 
$g_m\to f$ in $\cf$ as $m\to \infty$.
Here, note that the constant $c$ is taken independently of $m$, which derives from the fact that $c$ depends only on $c'$ in the proof of Lemma~\ref{lem:harmext}.
We first obtain that $g_m$ converges weakly to $f$ in $\cf$ and $\limsup_{m\to\infty}\ce(g_m-f)=\limsup_{m\to\infty}\ce(g_m)-\ce(f)\le0$.
Therefore, $\ce(g_m-f)\to0$ as $m\to\infty$.
By (B1), $g_m-f-\int_K(g_m-f)\,d\mu$ converges to $0$ in $L^2(K)$.
Since $\|g_m-f\|_{L^2(K)}\le c\|f\|_\cf+\|f\|_{L^2(K)}$, we have $\int_K(g_m-f)\,d\mu\to0$ as $m\to\infty$, which implies that $\|g_m-f\|_{L^2(K)}\to0$ as $m\to\infty$.
Thus, $g_m\to f$ in $\cf$ as $m\to\infty$.

Let $f_m=g_m-g_{m-1}$ where we set $g_{-1}\equiv 0$. Then, $f=\sum_{m=0}^\infty f_m$.
Since $f_i\in\cf(I^{b_0 j+b},\Phi(\hat I^{b_0 j+b}))$ for $i>j$
and $f_j\in\ch(I^{b_0 j+b},\Phi(\hat I^{B_0 j+b}))$, we have
$\ce (f_i,f_j)=0$ for $i\ne j$, so that 
\begin{equation}\label{eq:orthd}
\ce (f)=\sum_{m=0}^\infty \ce (f_m).
\end{equation}
Now, for each $f\in \cf$,
\begin{eqnarray}
(E_{(b_0i+b)}(Q_{b_0i+b}(\tilde f|_L)))^{1/2}&=& (E_{(b_0i+b)}(Q_{b_0i+b}(\tilde g_i|_L)))^{1/2}=
\left(E_{(b_0i+b)}\left(\sum_{j=0}^i Q_{b_0i+b}(\tilde f_j|_L)\right)\right)^{1/2}\nonumber\\
&\le  & \sum_{j=0}^i(E_{(b_0i+b)}(Q_{b_0i+b}(\tilde f_j|_L)))^{1/2}
\le  \sum_{j=0}^i(c_1\rho^{-b_0i-b}\ce_{\Phi(\hat I^{b_0i+b})}(f_j))^{1/2}\nonumber\\
&\le & \sum_{j=0}^i(c_1\rho^{-b_0i-b}c_0^{i-j}\ce_{\Phi(\hat I^{b_0j+b})}(f_j))^{1/2}\nonumber\\
&\le& \sum_{j=0}^i(c_1\rho^{-b_0i-b}c_0^{i-j}\ce (f_j))^{1/2},\label{eq:2.2}
\end{eqnarray}
where we apply Minkowski's inequality in the first inequality, 
(\ref{eq:qnapp}) in the second inequality, and
 Proposition \ref{theo:keyhin} in the
third inequality.

Applying (\ref{eq:2.2}) and noting that $\alpha^{d_w-d_f}=\rho$, we have 
\begin{eqnarray*}
\sum_{i=0}^\infty \alpha^{(d_w-d_f)(b_0i+b)}E_{(b_0i+b)}(Q_{b_0i+b}(\tilde f|_L))
&\le &\sum_{i=0}^\infty\rho^{b_0i+b}
\left(\sum_{j=0}^i (c_1\rho^{-b_0i-b}c_0^{i-j}\ce (f_j))^{1/2}\right)^2\\
&=&c_1\sum_{i=0}^\infty c_0^i\left(\sum_{j=0}^i (c_0^{-j}\ce (f_j))^{1/2}\right)^2\\
&\le &c_2\sum_{j=0}^\infty c_0^jc_0^{-j}\ce (f_j)= 
c_2\sum_{j=0}^\infty\ce (f_j)=c_2\ce (f).
\end{eqnarray*}
Here we used (\ref{eq:harheld})
in the second inequality and (\ref{eq:orthd}) 
in the last equality. Thus, we have
\[\sum_{n=0}^\infty \alpha^{(d_w-d_f)n}E_{(n)}(Q_{n}(\tilde f|_L))\le b_0c_2\ce(f).
\]
Combining this with Lemma \ref{theo:lem2} and (A8)', we have 
$\lll \tilde f|_L\lll_{\Lambda^\beta_{2,2}(L)}\le c_3\|f\|_\cf$, so that $\cf|_L\subset \Lambda_{2,2}^\beta(L)$
and $(\cf\cap C(K))|_L\subset \Lambda_{2,2}^\beta(L)\cap C(L)$. 
Noting that $\cf\cap C(K)$ is dense in $\cf$ due to the regularity of $(\ce,\cf)$, 
the claim follows by a simple limiting procedure.
\QED
\begin{rem} Even if (B4) does not hold, 
$\cf|_L\subset \hat\Lambda^\beta_{2,\infty}(L)$ hold.
Indeed, for each $f\in\cf$ and $n\in\bz_+$, we have  
by Lemma~\ref{lem:qnapp}, 
\begin{eqnarray*}
 c_1\ce(f)\ge c_1\ce_{\Phi(\hat I_n)}(f)
 \ge \rho^{n} E_{(n)}(Q_n(\tilde f|_L))
 = \alpha^{(2\beta-d)n} E_{(n)}(Q_n(\tilde f|_L)).
\end{eqnarray*}
Therefore, we have
\[
\left\|\left\{\alpha^{n\beta}\left(\alpha^{-nd}E_{(n)}(Q_n (\tilde f|_L))
\right)^{1/2}\right\}_{n=0}^\infty\right\|_{l^\infty}\le c_1^{1/2}\ce(f)^{1/2},
\]
so the same argument as above gives the result. 
\end{rem}
\subsection{Proof of the extension theorem}
In this section, we assume (A1)--(A8) and (C1)--(C2), and 
prove Theorem~\ref{theo:mainthm2}.
The conditions (C1) and (C2) will be defined below.

In order to construct an extension map $\xi$, we first define a Whitney-type decomposition and
an associated partition of unity. 
Let $\Omega^{(n)}=\bigcup_{m=0}^nI^m$ for $n\in\bz_+$.
For $w\in I^0=\{\emptyset\}$, set $A_w=W\setminus\cn_2(I)$ and $B_w=W\setminus \cn_1(I)$.
For $w\in I^n$ with $n\in\bn$,
set $A_w=(\cn_2(w)\cdot W)\setminus \cn_2(I^{n+1})$, $\hat A_w=\cn_2(w)\cdot W$, $B_w=(\cn_3(w)\cdot W)\setminus\cn_1(I^{n+1})$, and $\hat B_w=\cn_3(w)\cdot W$. 
Clearly
$K_{A_w}\subset K_{B_w}$, $K_{\hat A_w}\subset K_{\hat B_w}$, $K_{A_w}\cap K_{W^{|w|+1}\setminus B_w}=\emptyset$, and $K_{\hat A_w}\cap K_{W^{|w|+1}\setminus \hat B_w}=\emptyset$. 

By (A3), we see the following for $w,w'\in\bigcup_{n\in\bz_+}I^n$: 
\begin{eqnarray}
& c_1\alpha^{-|w|}\le \sd(L,K_{B_w})\le c_2\alpha^{-|w|}\mbox{ if }B_w\ne\emptyset,\label{eq:whi1}\\
&\mbox{there exists $l>0$ such that 
if $|w'|\ge|w|+l$ then $K_{B_w}\cap K_{\hat B_{w'}}= \emptyset$.}\label{eq:whi2} 
\end{eqnarray}
For $n\in\bn$ and $w\in\Omega^{(n)}$, we set
\[
  A^{(n)}_w=\left\{\begin{array}{ll}A_w&\mbox{if }|w|<n\\
  \hat A_w&\mbox{if }|w|=n\end{array}\right.,\quad
  B^{(n)}_w=\left\{\begin{array}{ll}B_w&\mbox{if }|w|<n\\
  \hat B_w&\mbox{if }|w|=n\end{array}\right.,
\]
and $R^{(n)}_w=\{w'\in\Omega^{(n)}: K_{B^{(n)}_w}\cap K_{B^{(n)}_{w'}}\ne \emptyset\}$.
We assume the following.
\begin{enumerate}
\item[(C1)] There exists a finite subset $\Gamma$ of $\bigcup_{n\in\bn}(\{n\}\times\Omega^{(n)})$ such that, for any $n\in\bn$ and $w\in\Omega^{(n)}$, there exist $(m,v)\in\Gamma$, a bijection $\iota\colon R^{(n)}_w\to R^{(m)}_v$, 
and a homeomorphism $F\colon K_{\bigcup_{u\in R^{(n)}_w}B^{(n)}_u}\to K_{\bigcup_{u\in R^{(m)}_v}B^{(m)}_u}$ satisfying that 
for every $u\in R^{(n)}_w$, $A^{(n)}_u$ and $A^{(m)}_{\iota(u)}$ are of the same type and so are $B^{(n)}_u$ and $B^{(m)}_{\iota(u)}$, for the homeomorphism $F$.
\end{enumerate}
For each $(m,v)\in \Gamma$, take a function $\bar\vp^{(m)}_v\in \cf\cap C(K)$ such that $0\le \bar\vp^{(m)}_v\le1$, $\bar\vp^{(m)}_v(x)=1$ on $K_{A^{(m)}_v}$, and $\bar\vp^{(m)}_v(x)=0$ on $K_{W^{|v|+1}\setminus B^{(m)}_v}$.
Such a function exists since $(\ce,\cf)$ is regular.
For $n\in\bn$ and $w\in\Omega^{(n)}$,  define  
$\vp^{(n)}_w(x)=\left\{\begin{array}{cl} \bar\vp^{(m)}_v( F(x))& \mbox {if } x\in B^{(n)}_w\\
0 &\mbox{otherwise}\end{array}\right.$,
where $m$, $v$ and $F$ are given in (C1).
We assume
\begin{enumerate}
\item[(C2)] $\vp^{(n)}_w\in\cf\cap C(K)$ for every $n\in\bn$ and $w\in\Omega^{(n)}$.
\end{enumerate}
For $n\in\bn$ and $w\in\Omega^{(n)}$, define
\[
\psi^{(n)}_w(x)=\frac{\vp^{(n)}_w(x)}{\sum_{w'\in\Omega^{(n)}}\vp^{(n)}_{w'}(x)},
\quad x\in K.
\]
This is well-defined since the sum in the denominator 
is not less than $1$.
$\psi^{(n)}_w$ is continuous and takes values between 0 and 1.
Since $\vp^{(n)}_w\in\cf$ and vanishes outside of $K_{B^{(n)}_w}$, so does $\psi^{(n)}_w$.
For each $f\in \Lambda^\beta_{2,2}(L)\cap C(L)$, define
\[
\xi^{(n)} f(x)=
\sum_{w\in\Omega^{(n)}}\psi^{(n)}_{w}(x)Q_{|w|}f(w)=\sum_{w\in\Omega^{(n)}}\psi^{(n)}_{w}(x)\mint_{L_w}f(s)d\nu(s).
\]
$\xi^{(n)}$ is a linear map from $\Lambda^\beta_{2,2}(L)\cap C(L)$ to $\cf\cap C(K)$.
For $x\in K\setminus L$, $\xi^{(n)}f(x)$ is independent of $n$ if $n$ is sufficiently large because of (\ref{eq:whi2}).
Therefore, for $f\in \Lambda^\beta_{2,2}(L)\cap C(L)$,
\begin{equation}\label{eq:xidef}
\xi f(x):=\left\{
\begin{array}{ll}\displaystyle\lim_{n\to\infty}\xi^{(n)}f(x),& x\in K\setminus L\\
 f(x),& x\in L\end{array}
\right.
\end{equation}
is well-defined and $\xi^{(n)}f$ converges to $\xi f$ $\mu$-a.e.
\Proofwithname{Proof of Theorem~\ref{theo:mainthm2}.}
We first prove that $\xi f$ is continuous on $K$.
Since  $\xi f$ is continuous on $K\setminus L$ by the construction, it is enough to show for each $x_0\in L$ that 
\begin{equation}
\lim_{\begin{subarray}{c}{x\to x_0}\\{x\in K\setminus L}\end{subarray}}\xi f(x)=f(x_0).
\label{eq:need1}\end{equation} 
Since $f$ is uniformly continuous on $L$, if we set 
$\omega_a(f)=\sup\{|f(s)-f(t)|:s,t\in L,\ \sd(s,t)\le a\}$ 
for $a>0$, then $\lim_{a\to0}\omega_a(f)=0$.
Let $x_0\in L$, $x\in K\setminus L$, and $\delta=\sd(x,x_0)$. 
Suppose that $w\in\bigcup_{n\in\bn}I^n$ satisfies $x\in K_{B_w}$. 
Then,
$c_1\alpha^{-|w|}\le \sd(L,K_{B_w})\le \sd(x_0,x)=\delta$ by $(\ref{eq:whi1})$. 
Next, take $y\in L_w$ and choose $z\in 
K_{B_w}$ that satisfies $\sd(y,z)=\sd(y,K_{B_w})\le c_2\alpha^{-|w|}$. Then, since 
$\mathop{\rm diam\,}(K_{B_w})\asymp \alpha^{-|w|}$, we have 
\[
\sd(y,x_0)\le  \sd(y,z)+\sd(z,x)+\sd(x,x_0)
\le  c_2\alpha^{-|w|}+c_3\alpha^{-|w|}+\delta\le c_4\delta.
\]
Therefore, $\mint_{L_w}|f(y)-f(x_0)|\,d\nu(y)\le \omega_{c_4 \delta}(f)$.
Now, take $n$ sufficiently large so that $x\notin\bigcup_{w\in I^n}K_{\hat B_w}$.
Then, $\xi^{(n)}f(x)=\xi f(x)$ and
\begin{eqnarray*}
|\xi f(x)-f(x_0)|
&=& \left|\sum_{w\in\Omega^{(n)}} \psi^{(n)}_{w}(x) 
\mint_{L_w}(f(y)-f(x_0))\,d\nu (y)\right|\\
&\le&\sum_{w\in\Omega^{(n-1)},\ x\in K_{B_w}}\psi^{(n)}_{w}(x) 
\mint_{L_w}|f(y)-f(x_0)|\,d\nu (y)\\
&\le& \omega_{c_4\delta}(f).
\end{eqnarray*}
Thus $(\ref{eq:need1})$ is proved.

Next, we will prove $\{\xi^{(n)}f\}_{n\in\bn}$ is bounded in $\cf$.
Noting that $\int_K \psi^{(n)}_w(x)\,d\mu(x)\le c_5 \alpha^{-d_f|w|}$ for all $n\in\bn$ and $w\in\Omega^{(n)}$ for some $c_5>0$, we have
\begin{eqnarray}
\|\xi^{(n)} f\|_{L^2(K,\mu)}^2
&=&\int_K\left(\sum_{w\in \Omega^{(n)}}\psi^{(n)}_{w}(x)\mint_{L_w}
f(s)\,d\nu (s)\right)^2d\mu (x)\nonumber\\
&\le &\int_K\left(\sum_{w\in \Omega^{(n)}}\psi^{(n)}_{w}(x)\mint_{L_w}
f(s)^2\,d\nu (s)\right)d\mu (x)\nonumber\\
&\le &\sum_{w\in \Omega^{(n)}}c_5 \alpha^{-d_f |w|}\alpha^{d |w|}\int_{L_w}f(s)^2d\nu (s)\nonumber\\
&=& c_5 \sum_{k=0}^n \alpha^{(d-d_f)k}\|f\|^2_{L^2(L,\nu)}\nonumber\\
&\le& \frac{c_5}{1-\alpha^{d-d_f}}\|f\|^2_{L^2(L,\nu)}.
\label{eq:l2hk}\end{eqnarray}

For $n\in\bn$, $w\in\Omega^{(n)}$ with $m=|w|$, 
let $\bar R^{(n)}_w=\bigcup_{v\in R^{(n)}_w} v\cdot I^{m+l-|v|}\subset I^{m+l}$, where $l$ is provided in (\ref{eq:whi2}).
For $g\in L^2(L,\nu)$, we define
\[
E^{(n)}_w(g)=\sum_{\begin{subarray}{c}u,v\in \bar R^{(n)}_w\\u\stackrel{m+l,L}{\longleftrightarrow}v\end{subarray}}(Q_{m+l} g(u)-Q_{m+l}g(v))^2,\quad
\bar E^{(n)}_w(g)= \ce_{\Phi(B^{(n)}_w)}\left(\sum_{v\in R^{(n)}_w}Q_{|v|}g(v)\psi^{(n)}_v\right).
\]
Both $E^{(n)}_w(g)$ and $\bar E^{(n)}_w(g)$ are determined only by the values $\{Q_{m+l}g(u)\}_{u\in \bar R^{(n)}_w}$.
If $E^{(n)}_w(g)=0$, then $Q_{m+l}g$ is constant on $\bar R^{(n)}_w$, which implies that $\bar E^{(n)}_w(g)=0$.
Therefore, there exists $c^{(n)}_w>0$
such that $\bar E^{(n)}_w(g)\le c^{(n)}_w E^{(n)}_w(g)$ for every $g\in\cf$.
Due to (C1) and Lemma~\ref{lem:scale}, there exists some $c_6>0$ such that 
$\bar E^{(n)}_w(g)\le c_6 \rho^{|w|} E^{(n)}_w(g)$ for all $n\in\bn$, $w\in\Omega^{(n)}$ and $g\in\cf$.
It also holds that there exists $c_7>0$ 
independent of $m$ such that 
$\sum_{w\in I^m}E^{(n)}_w(g)\le c_7 E_{(m+l)}(Q_{m+l}g)$ for 
all $n$ and $g\in L^2(L,\nu)$.
Then, we have
\begin{eqnarray*}
\ce(\xi^{(n)}f)
&\le& \sum_{w\in\Omega^{(n)}}\ce_{\Phi(B^{(n)}_w)}(\xi^{(n)}f)
= \sum_{w\in\Omega^{(n)}}\bar E^{(n)}_w(f)\\
&\le& c_6\sum_{m=0}^n\sum_{w\in I^m}\rho^m E^{(n)}_w(f)
\le c_6 c_7\sum_{m=0}^n\rho^m E_{(m+l)}(Q_{m+l}f)\\
&\le& c_8\sum_{m=0}^{\infty}\rho^m E_{(m)}(Q_mf).
\end{eqnarray*}
Since $\alpha^{2\beta-d}=\alpha^{d_w-d_f}=\rho$, we obtain $\ce(\xi^{(n)}f)\le c_8 \|f\|_{\Lambda^\beta_{2,2}(L)}^2$ by Lemma \ref{theo:lem2}.

By combining this with (\ref{eq:l2hk}), $\{\xi^{(n)}f\}_{n\in\bn}$ is bounded in $\cf$ and we conclude that $\xi f\in\cf$ and 
$\|\xi f\|_\cf\le c_9\|f\|_{\Lambda^\beta_{2,2}(L)}$ for some $c_9>0$.

Next,
take any $\Lambda^\beta_{2,2}(L)$-Cauchy sequence 
$\{f_n\}_{n\in\bn}\subset \Lambda^\beta_{2,2}(L)\cap C(L)$ and let $f\in  \Lambda^\beta_{2,2}(L)$
be the limit point. By the above result, $\{\xi f_n\}_{n\in\bn}\subset \cf\cap C(K)$
is a $\ce_1$-Cauchy sequence. Let $g\in \cf$ be the limit point. 
Since $\xi f_n|_L=f_n$ 
and a subsequence $\xi f_{n_k}$ converges to $\tilde g$ q.e.,
$\tilde g|_L=f$ $\nu$-a.e. Thus, $\xi$ can extend to a continuous map from  $\hat\Lambda^\beta_{2,2}(L)$ to $\cf$ such that $\widetilde{\xi f}|_L=f$ $\nu$-a.e.\ for 
$f\in \hat\Lambda^\beta_{2,2}(L)$. 
\QED
\begin{rem}\label{thm:finunil}
Let $\{L_i\}_{i=1}^m$ be a finite number of self-similar subsets of $K$ where each $L_i$ 
is constructed by the same number of contraction maps and satisfies (A2), 
the second identity of (A4),
(A7), and (A8) in Section~2. Let $L=\bigcup_{i=1}^mL_i$. 
With suitable changes for $A_w$, $B_w$ etc., we 
can consider conditions
(C1)$^*$--(C2)$^*$ as the corresponding 
(C1)--(C2). 
Define $\Lambda^\beta_{2,2}(L)$ as in Definition~\ref{defn:Besov}.
Then, under such conditions, 
Theorem~\ref{theo:mainthm2} 
is still valid, i.e.\ there is a linear map $\xi$ from $\Lambda^\beta_{2,2}(L)$ 
to $\cf$ such that $\xi\bigl(\Lambda^\beta_{2,2}(L)
\cap C(L)\bigr)\subset \cf\cap C(K)$, $\widetilde{\xi f}|_L=f$ and
\[
\|\xi f\|_\cf\le c_1\sum_{i=1}^m\lll f|_{L_i}\lll_{\Lambda^\beta_{2,2}(L_i)},
\quad f\in\Lambda^\beta_{2,2}(L).
\]
\end{rem}
\section{Complementary results}
In this section, we give sufficient conditions concerning (A8) and (B3),
and discuss a suitable choice of $\cf_A$ for $A\subset W^m$. 
We first define fractional diffusions in the sense of \cite{bar} Definition 3.5.
\begin{defn}\label{fradif}
Let $(X,\sd)$ be a complete metric space where $\sd$ has the midpoint property; for each
$x,y\in X$, there exists $z\in X$ such that $\sd(x,y)=\sd(x,z)/2=\sd(z,y)/2$. 
For simplicity, we assume $\mbox{diam}\,X=1$.
Let $\mu$ be a Borel measure on $X$ such that there exists $d_f>0$ with $\mu (B(x,r))\asymp r^{d_f}$
for all $0<r\le 1$. A Markov process $\{Y_t\}_{t\ge 0}$ is a fractional diffusion on $X$ if\\
1) $Y$ is a $\mu$-symmetric conservative Feller diffusion,\\
2) $Y$ has a symmetric jointly continuous transition density 
$p_t(x,y)~(t>0, x,y\in X)$ which satisfies the Chapman-Kolmogorov equations and 
has the following estimate,
\begin{eqnarray*}
&&c_{1}t^{-d_f/d_w}\exp (-c_{2}(\sd(x,y)^{d_w}t^{-1})^{1/(d_w-1)})
\le p_t(x,y)\\
&\le & c_{3}t^{-d_f/d_w}\exp (-c_{4}(\sd(x,y)^{d_w}t^{-1})^{1/(d_w-1)})
~~\qquad\qquad\mbox{ for all}~~0<t<1,~x,y\in X,
\end{eqnarray*}
with some constant $d_w\ge 2$. 
\end{defn}

\begin{pro}\label{thm:pr4.2} 
\rom{(A8)} holds for the following 
three cases.\\ 
\rom{1)} There exists $c>0$ such that $\|f\|_{L^\infty(K)}\le c\|f\|_\cf$ for all $f\in\cf$.\\
\rom{2)} The diffusion process corresponding to $(\ce,\cf)$
is the fractional diffusion and \rom{(A7)} holds.\\
\rom{3)} $K\subset \br^n$, 
\rom{(A7)} holds, and $\cf=\Lambda^{d_w/2}_{2,\infty}(K)$.
\end{pro}
\Proof
Suppose that 1) holds.
Then, for any nonempty set $D$ of $K$, $\Cp(D)\ge c^{-2}$.
Therefore, $\nu(D)\le 1\le c^2\Cp(D)$.

The proof when 2) holds is 
similar to Lemma 2.5 of \cite{bk}, 
but we will give it for completeness. 
Let $g_1(\cdot,\cdot)$ be the $1$-order 
Green density given by  
\[E^x\left[\int_0^{\infty}e^{- t}f(X_t)dt\right]
=\int_K  g_1(x,y)f(y)d\mu,\]
for any Borel measurable function $f$, 
where $\{X_t\}$ is the diffusion corresponding to $(\ce, \cf)$. Then, 
since $\{X_t\}$ is a fractional diffusion, we have 
\begin{eqnarray}
g_1(x,y)&\asymp &  \left\{ \begin{array}{ll} c_{1}\sd(x,y)^{d_w-d_f}
&~~\mbox{if}~d_f> d_w,\label{eq:green}\\ 
-c_{2}\log \sd(x,y)+c_{3} &~~\mbox{if}~d_f=d_w,
\\ c_{4} &~~\mbox{if}~d_f<d_w.
\end{array}
\right.\end{eqnarray}
See \cite{bar} Proposition 3.28 for the proof. 
If $d_f<d_w$ then points have strictly positive capacity, and the result is
immediate. We prove the result for $d_f>d_w$: the proof for $d_f=d_w$
is similar.
It is well-known that for each compact set $M\subset K$,
\begin{equation}
\Cp(M)=\sup \left\{m (M):\begin{array}{l} m\mbox{ is a positive Radon measure, }
\supp m\subset M,\\
G_1 m(x)\equiv \int_M g_1(x,y) m(dy)\le 1\mbox{ for every } x\in K
\end{array}
\right\}.
\label{eq:capca}
\end{equation}
Using the above estimates of $g_1(\cdot,\cdot)$,
\begin{eqnarray*}
\int_M  g_1(x,y) \nu (dy)  &\le &
\int_K g_1(x,y) \nu (dy)
\le \sum_{n=0}^{\infty}\int_{\alpha^{-n-1}
  \le \sd(x,y) <\alpha^{-n}} g_1(x,y)\nu (dy)\\
&\le &c_5\sum_n \alpha^{n(d_f-d_w)}\nu ( \alpha^{-n-1}\le \sd(x,y) <\alpha^{-n} )
\le  c_6\sum_{n} \alpha^{n(d_f-d_w-d)} \equiv c_7<\infty,
\end{eqnarray*}
because of the assumption $d_f-d<d_w$. 
Thus, setting $\nu_M (\cdot)\equiv \nu (\cdot \cap M)$,
we have $G_1 \nu_M\le c_7$. Using $(\ref{eq:capca})$, $\mbox {Cap}(M)\ge \nu
(M)/c_7$ for each compact set $M$.

For 
3), we will use the results by Jonsson-Wallin in \cite{jw} and by Triebel in 
\cite{tr}. Denote the Lipschitz and the
Besov spaces in the sense of Jonsson-Wallin by
$\mbox{Lip}_{JW} (\alpha, p,q, K)$ and $B_{\alpha, JW}^{p,q}(K)$
(see page 122--123 in \cite{jw} for definition). Note that 
$\mbox{Lip}_{JW} (\alpha, p,q, K)\subset B_{\alpha, JW}^{p,q}(K)$
and they are equal when $\alpha\notin \bn$ (page 125 in \cite{jw}). 
For each $f\in \Lambda^{d_w/2}_{2,\infty}(K)$,
$(f,0,\ldots,0)\in \mbox{Lip}_{JW} (d_w/2, 2,\infty, K)$. Thus, using 
the extension theorem in page 155 of \cite{jw}, we have
\[ \Lambda^{d_w/2}_{2,\infty}(K)\subset \mbox{Lip}_{JW} (d_w/2, 2,\infty, K)
\subset B_{d_w/2, JW}^{2,2}(K)\subset \Lambda_\gamma^{2,\infty}(\br^n)|_K,
\]
where $\gamma=(d_w+n-d_f)/2$ and $\Lambda_\gamma^{p,q}(\br^n)$ is a classical Besov space
on $\br^n$. Now, since $d_w-d_f>-d$ (due to (A7)), $\Lambda_\gamma^{2,\infty}(\br^n)
\subset \Lambda_{(n-d)/2}^{2,1}(\br^n)$. Finally, by Corollary 18.12 (i) in \cite{tr},
we have $\mbox{tr}_L\;\Lambda_{(n-d)/2}^{2,1}(\br^n)=L^2(L,\nu)$. (Note that this trace 
in the sense of Triebel is simply restriction and there is no corresponding extension.)
Combining these facts, we have $ \cf|_L\subset L^2(L,\nu)$, which means 
$\|\tilde f|_L\|_{L^2(L,\nu)}\le c_9\|f\|_\cf$ for all $f\in \cf$. 
Therefore, (A8)' holds.
\QED
We now make one concrete choice of $\cf_A$ for $A\subset W^m$ 
and show that such a choice is suitable for Dirichlet forms whose corresponding 
processes are the fractional diffusions. 
By (A5), (A6), and the self-similarity of $\mu$, for any $w\in\bigcup_{n\in\bz_+}W^n$, there exists $c>0$ such that $\Cp(D)\le c\Cp(F_w(D))$ for any $D\subset K$.
We assume the converse as follows.
\begin{enumerate}
\item[(A*)]
For any $w\in\bigcup_{n\in\bz_+}W^n$, there exists $c>0$ such that $\Cp(F_w(D))\le c\Cp(D)$ for any $D\subset K$.
\end{enumerate}
For a subset $A$ of $W^m$ for some $m\in\bn$, 
we say that 
a collection $\{f_w\}_{w\in A}$ of functions in $\cf$
is compatible if $\tilde f_v(F_v^{-1}(x))=\tilde f_w(F_w^{-1}(x))$ q.e.\ on $K_v\cap K_w$ for every $v,w\in A$. 
Note that this is well-defined by (A*). Define
\begin{equation}\label{eq:cfadef}
  \cf_A=\{f\in L^2(K_A,\mu|_{K_A}): F_w^* f\in\cf \mbox{ for all }w\in A \mbox{ and $\{F_w^* f\}_{w\in A}$ is compatible}\}.
\end{equation}
If we equip $A$ with a graph structure so that $v\in A$ and $w\in A$ are connected if $\Cp(K_v\cap K_w)>0$, then $A$ is $\ce$-connected when $A$ is a connected graph.
This is verified by using (B1).
\begin{lem}\label{lem:good}
  For $A\subset W^m$ with $m\in\bz_+$,
  $(\ce_A,\cf_A)$ is a strong local Dirichlet form on 
  $L^2(K_A,\mu|_{K_A})$.
\end{lem}
\Proof
  Let $\{f_n\}_{n\in\bn}\subset \cf_A$ be a Cauchy sequence in $\cf_A$.
  Let $g$ be the limit in $L^2(K_A)$.
  Let $w\in A$. Since $\{F_w^* f_n\}_{n\in\bn}$ is a Cauchy sequence in $\cf$, $F_w^* f_n\to\ F_w^* g$ in $\cf$.
  It is also easily deduced that 
  $\{F_w^* g\}_{w\in A}$ 
is compatible. Therefore, $g\in \cf_A$ and $f_n\to g$ in $\cf_A$.
  This implies that $(\ce_A,\cf_A)$ is a closed form on $L^2(K_A)$.
  The Markov property and the strong locality are inherited from those of $(\ce,\cf)$ via relation (\ref{eq:energy}).
\QED

\begin{cor}\label{cor:a*co}
 Assume case \rom{1)} or \rom{2)} in Proposition~\ref{thm:pr4.2}.
 Then, \rom{(A*)} holds. 
\end{cor}
\Proof
In case~1), non-empty sets have uniform positive capacities, which implies (A*).
In case~2), (A*) is an easy consequence of (\ref{eq:green}) and (\ref{eq:capca}).
\QED

In the rest of this section, we will discuss 
sufficient conditions for \rom{(B3)}. 

\begin{lem}\label{lem:Dw}
Let $A\subset W_m$, $m\in\bz_+$.
Then, $\cf_A$ is compactly imbedded in 
$L^2(K_A,\mu|_{K_A})$.
Suppose that $A$ is $\ce_A$-connected.
Then, when we set $\ca=\{f\in\cf_A: \int_{K_A}f\,d\mu=0,\ \ce_A(f)\le C\}$ for a constant $C>0$, $\ca$ is bounded in $\cf_A$.
\end{lem}
\Proof
  Let $\cb$ be a bounded subset of $\cf_A$.
  For each $v\in A$, $\{F^*_v f: f\in\cb\}$ is bounded in $\cf$.
  By (B1), we can take a sequence $\{f_n\}_{n\in\bn}$ from $\cb$ such that $F^*_v f_n$ converges in 
  $L^2(K)$.
  Therefore, we can take a sequence from $\cb$ converging in $L^2(K_A)$. 
  This implies the first assertion.
  
  By combining this with the $\ce_A$-connectedness of $K_A$, there exists $c>0$ such that $\|f-{}-\nobreak\!\!\!\!\!
  \int_{K_A}f\,d\mu\|_{L^2(K_A)}^2 \le c\ce_A(f)$ for every $f\in\cf_A$.
  The latter assertion follows from this immediately.
\QED
We now give a sufficient condition for \rom{(B3)~(2)}.
\begin{pro}\label{prop:ehi}
  The following condition \rom{(EHI1)} implies \rom{(B3)~(2)}.
  \begin{enumerate}
  \item[\rom{(EHI1)}] For any $v\in \Xi$, there exist some $c_1>0$,
  subsets $D''(v)$ and $D'''(v)$ of $D^\sharp(v)$ such that
  $D'''(v)\subset D''(v)\subset D^\sharp(v)$, 
  $K_{D''(v)}\cap K_{S^{|v|}\setminus D^\sharp(v)}=\emptyset$,
  $K_v\cap K_{S^{|v|}\setminus D'''(v)}=\emptyset$ and
  $\esssup_{x\in K_{D'''(v)}} h(x) \le c_1\essinf_{x\in K_{D'''(v)}} h(x)$ for every $h\in \ch(D^\sharp(v))$ with $h\ge0$ $\mu$-a.e.
  \end{enumerate}
\end{pro}
\Proof
First, we apply Lemma~\ref{lem:harmext} to $g\in\cf$ with $A=D^\sharp(v)$ and 
$J=\emptyset$, and denote $g'$ there by $Hg$.
We will follow the proof of Theorem~2.2 in \cite{hin}. 
For $h\in \ch(D^\sharp(v))$ with $h\ge0$ $\mu$-a.e., we have, by (EHI1),
\[
  \esssup_{x\in K_{D'''(v)}} h(x) \le c_1\essinf_{x\in K_{D'''(v)}} h(x)
  \le c_2\| h\|_{L^2(K_{D'''(v)})}.
\]
For $h\in \ch(D^\sharp(v))$, let $h_+(x)=\max\{h(x),0\}$ and $h_-(x)=\max\{-h(x),0\}$.
Since $h=Hh=Hh_+-Hh_-$ and $Hh_\pm\ge0$ $\mu$-a.e., we have
\begin{eqnarray}
\esssup_{x\in K_{D'''(v)}} |h(x)|
&\le&\esssup_{x\in K_{D'''(v)}} Hh_+(x)+\esssup_{x\in K_{D'''(v)}} Hh_-(x)\nonumber\\
&\le& c_2(\|Hh_+\|_{L^2(K_{D'''(v)})}+\|Hh_-\|_{L^2(K_{D'''(v)})})\nonumber\\
&\le& c_3(\|h_+\|_{\cf_{D''(v)}}+\|h_-\|_{\cf_{D''(v)}})\nonumber\\
&\le&2c_3\|h\|_{\cf_{D''(v)}}.
\label{eq:dominationw}
\end{eqnarray}
In order to prove (B3) (2), it suffices to prove the following.
\begin{itemize}
\item[$(*)$] If a sequence $\{h_l\}$ in $\ch(D^\sharp(v))$ converges weakly to 0 in $\cf_{D^\sharp(v)}$, then there exists a subsequence $\{h_{l(k)}\}$ such that $F_v^* h_{l(k)}$ converges strongly to $0$ in $\cf$.
\end{itemize}
Indeed, suppose $(*)$ holds. 
Let $\{f_m\}$ be a sequence in $\ch(D^\sharp(v))$ that is bounded in $\cf_{D^\sharp(v)}$.
We can take a subsequence $\{f_{m(l)}\}$ and $f\in\cf_{D^\sharp(v)}$ such that $f_{m(l)}$ converges weakly to $f$ in $\cf_{D^\sharp(v)}$. 
Take $g_l\in\ch(D^\sharp(v))$ such that $g_l\to f$ in $\cf_{D^\sharp(v)}$.
Applying $(*)$ to $h_l:=f_{m(l)}-g_l$, we can take a sequence $\{l(k)\}$ diverging to $\infty$ such that 
$F^*_v f_{m(l(k))}\to F^*_v f$ in $\cf$. 
This implies (B3)~(2).

In order to prove $(*)$, recall the notion of the energy measure. For $f\in \cf\cap L^\infty(K)$, the energy measure $\mu_{\langle f\rangle}$ is a unique positive Radon measure on $K$ such that the following identity holds for every $g\in\cf\cap C(K)$:
\[
  \int_{K} g\,d\mu_{\langle f\rangle}=2\ce(f,fg)-\ce(f^2,g).
\]
Now, by (\ref{eq:dominationw}), $C:=\esssup_{x\in K_{D'''(v)}} |h_l(x)|$ is bounded in $l$.
Define $\hat h_{l}=(-C)\vee h_l \wedge C$.
Since $\cf_{D^\sharp(v)}$ 
is compactly imbedded in $L^2(K_{D^\sharp(v)})$ by Lemma~\ref{lem:Dw}, $\{h_l\}$ converges to $0$ in $L^2(K_{D^\sharp(v)})$. 
Take a subsequence $\{h_{l'}\}$ converging to $0$ $\mu$-a.e.\ on $K_{D^\sharp(v)}$.
Since $(\ce,\cf)$ is regular, we can take $\vp\in\cf\cap C(K)$ such that $0\le\vp\le 1$ on $K$, 
$\vp=1$ on $K_w$, and $\vp=0$ outside $K_{D'''(v)}$.
We have
\[
  0=2\ce(h_{l'},\hat h_{l'}\vp)
  =2\ce(\hat h_{l'},\hat h_{l'}\vp)
  =\ce(\hat h_{l'}^2,\vp)+\int_K\vp\,d\mu_{\langle \hat h_{l'}\rangle},
\]
because $\hat h_{l'}\vp$ vanishes outside $K_{D'''(v)}$.
Note that $\ce(\hat h_{l'}^2)\le 4C^2 \ce(h_{l'})$, which is
bounded in $l'$. A suitable subsequence $\hat h_{l''}$ can be taken so that $\{\hat h_{l''}^2\}$ converges weakly to some $g$ in
$\cf$. 
Since $g=0$ on $K_{D^\sharp(v)}$, $\ce(\hat h_{l''}^2,\vp)\to\ce(g,\vp)=0$ as $l''\to\infty$.
On the other hand,
\begin{eqnarray*}
\int_{K}\vp\,d\mu_{\langle \hat h_{l''}\rangle}
&=&\sum_{z\in S^{|v|}} \rho^{|v|} \int_K F^*_z \vp\,d\mu_{\langle F^*_z \hat h_{l''}\rangle}\\
&\ge& \rho^{|v|} \int_K F^*_v \vp\,d\mu_{\langle F^*_v \hat h_{l''}\rangle}\\
&=& \rho^{|v|} \mu_{\langle F^*_v \hat h_{l''}\rangle}(K)
= 2\rho^{|v|} \ce(F^*_v \hat h_{l''})
= 2\rho^{|v|} \ce(F^*_v h_{l''}).
\end{eqnarray*}
Combining these estimates, we obtain 
$\varlimsup_{l''\to\infty}\ce(F^*_v h_{l''})\le 0$.
Therefore, $F^*_v h_{l''}$ converges to $0$ in $\cf$.
This proves $(*)$.
\QED
We next give sufficient conditions for \rom{(B3) (1) (b)}.
\begin{pro}\label{prop:wv}
  The following conditions imply \rom{(B3) (1) (b)}.
  \begin{enumerate}
  \item[$(1)$] $\cf=\Lambda^\beta_{2,\infty}(K)$ for some $\beta>0$.
  \item[$(2)$] For each $v\in \Xi$, $D'(v)$ is $\ce_{D'(v)}$-connected.  
  \item[$(3)$] For each $w\in\bigcup_{n\in\bz_+}\hat I^{n+m_0}$, there exist 
  subsets $D^\sharp(w)$, $D^{(1)}(w)$, $D^{(2)}(w)$ of $D'(w)$ such that $D^\sharp(w)\subset D^{(1)}(w)\subset D^{(2)}(w)\subset D'(w)$ and the following hold.
  \begin{enumerate}
  \item[\rom{(a)}] There exists $v\in \Xi$ such that both $D'(w)$ and $D'(v)$, and $D^\sharp(w)$ and $D^\sharp(v)$, are of the same type by the same map $F$.
  \item[\rom{(b)}] $K_{D^\sharp(w)}\cap K_{S^{|w|}\setminus D^{(1)}(w)}=K_{D^{(2)}(w)}\cap K_{S^{|w|}\setminus D'(w)}=\emptyset$.
  \item[\rom{(EHI2)}] There exists $c>0$ such that $\esssup_{x\in K_{D^{(1)}(w)}} h(x) \le c\essinf_{x\in K_{D^{(1)}(w)}} h(x)$ for $h\in \ch(D^{(2)}(w))$ with $h\ge0$ $\mu$-a.e.
 \end{enumerate}
 \end{enumerate}
\end{pro}
\Proof
  Let $g$ be a function in $\cf$ such that $\int_{K_{D'(w)}}g\,d\mu=0$ and $\rho^{|w|}\ce_{D'(w)}(g)\le 1$. Let $f(x)=g(F^{-1}(x))$, $x\in K_{D'(v)}$.
  Then, $f\in \cf_{D'(v)}$, $\int_{K_{D'(v)}}f\,d\mu=0$ and $\rho^{|v|}\ce_{D'(v)}(f)\le 1$.
  By Lemma~\ref{lem:Dw}, $\|f\|_{\cf(D'(v))}\le C$, where $C$ is a constant independent of $w$.
  Suppose moreover that $g\in\ch(I^{|w|},D'(w))$.
  Apply Lemma~\ref{lem:harmext} to $g$ with $A=D^{(2)}(w)$ and $J=\emptyset$ and denote $g'$ there by $g_1$.
  Let $g_2=g-g_1$.
  By (EHI2) and the same argument in the first part of the proof of Proposition~\ref{prop:ehi}, $g_1$ is bounded on $D^{(1)}(w)$.
  Take a function $\psi\in\cf$ such that $0\le\psi\le 1$, $\psi=0$ on $K_{S^{|w|}\setminus D^{(1)}(w)}$ and $\psi=1$ on $K_{D^\sharp(w)}$.
  Then, $g_1\psi\in\cf$.
  Since both $g_1\psi$ and $g_2$ vanish on $K_{S^{|w|}\setminus D^{(2)}(w)}$, when we set
  $f'(x)=\left\{\begin{array}{ll}(g_1\psi+g_2)(F^{-1}(x)),&x\in K_{D'(v)}\\0,&x\in K\setminus K_{D'(v)}\end{array}\right.$, $f'$ belongs to $\cf$ by using the fact $\cf=\Lambda^\beta_{2,\infty}(K)$.
  Since $f'=f$ on $K_{D^\sharp(v)}$, we have $f'\in\ch(I^{|v|},D^\sharp(v))$ and $\|f'\|_{\cf_v}=\|f\|_{\cf_v}\le C$.
  These conclude the assertion.
\QED
\section{Examples}
In this section, we choose $\cf_A$ as in (\ref{eq:cfadef}) for $A\subset W^m$.\\

\noindent
1) Sierpinski gaskets: Let $\{a_0,a_1,\ldots,a_n\}\subset \br^n$ be the vertices 
of $n$-dimensional simplex. 
Let $W=S=\{0,1,\ldots,n\}$ and let $F_i(x)=(x-a_i)/2+a_i$ for $x\in \br^n$ and 
$i=0,1,\ldots ,n$. 
Then the unique non-void compact set $K$ which satisfies $K=\bigcup_{i=0}^n F_i(K)$ 
is the $n$-dimensional Sierpinski gasket. 
The map $\Phi$ in Lemma~\ref{lem:Phi} is the identity map.
It is well-known (see \cite{bar,BP,kig} etc.) 
that there is a self-similar Dirichlet form $(\ce,\cf)$ on $L^2(K,\mu)$ 
where the corresponding diffusion is the fractional diffusion. 
In particular, $\cf=\Lambda_{2,\infty}^{d_w/2} (K)\subset C(K)$ where $d_w=(\log (n+3))/(\log 2)$.  
Note that $d_w-d_f>0$ in this case. 
Let $L$ be the $(n-1)$-dimensional gasket determined by $\{F_i\}_{i=0}^{n-1}$.
That is, $I=\{0,1,\ldots,n-1\}$. Let $\hat I=I$, $M=1$. 
It is easy to see that 
$\cf_A=\cf|_{K_A}$ for each $A\subset W^m$. 
Then (A1)--(A7), (B2), and (C1)--(C2)
are easy to check with $\rho=(n+3)/(n+1)$. (A8) holds by Proposition \ref{thm:pr4.2} and
(B1) holds by \cite{kig} Lemma 3.4.5. 
For (B3), define $l_0=0$, $m_0=0$, $D'(w)=\{w\}$ for $w\in\bigcup_{n\in\bz_+}I^n$, $\Xi=\{\emptyset\}$, and $D^\sharp(\emptyset)=\{\emptyset\}$.
It is easy to check (B3)~(1) by using Lemma~\ref{lem:scale} and Lemma~\ref{lem:Dw}.
Since $\ch_{D'(w)}(D'(w))$ is a finite dimensional space, 
(B3)~(2) is clearly true.
We will prove (B4).
Let $f\in\cf$ and $\ce_{S^m\setminus I^m}(f)=0$ for some $m\in\bz_+$.
Then, for each $w\in S^m\setminus I^m$,
$\ce(F^*_w f)=\rho^{-m}\ce_w(f)=0$.
Therefore, $f$ is constant on $K_w$ for each $w\in S^m\setminus I^m$.
We consider an unoriented graph with a vertex set $V=S^m\setminus I^m$ and an edge set $\{(v,w)\in V\times V: \Cp(K_v\cap K_w)>0\}$.
Then, $V$ is a connected set. Note that $(v,w)$ is an edge if and only if $K_v\cap K_w\ne\emptyset$.
Therefore, $f$ should be constant on $K_{S^m\setminus I^m}$, thus constant on $K\setminus L$.
This concludes that (B4) holds. 
Therefore, we have by Theorem~\ref{theo:mainthm1}, Theorem~\ref{theo:mainthm2} and Remark~\ref{rem:2.5}, 
\[\cf|_L=\Lambda_{2,2}^\beta (L)\qquad\mbox{where }~~\beta=\frac{d_w}2-
\frac{\log (1+1/n)}{2\log 2}.\]
When $n=2$, this relation was obtained in \cite{jo}.

\medskip

\noindent 2) Pentakun: 
Let $a_k=e^{2k\sqrt{-1}\pi/5+\sqrt{-1}\pi/2}\in\bc$, $k=0,1,2,3,4$.
Let $W=\{0,1,2,3,4,5,6\}$, $S=\{0,1,2,3,4\}$, $I=\hat I=\{2,3,5,6\}$ and $M=1$.
Let $\fg=\{G_k\}_{k=0}^4$ with $G_k:\bc\to\bc$ defined by $G_k(z)=e^{2k\sqrt{-1}\pi/5} z$.
For $i=0,1,2,3,4$, define a contraction map $F_i:\bc\to\bc$ by
$F_i(z)=\alpha^{-1}(z-a_i)+a_i$, where $\alpha=\frac{3+\sqrt5}{2}$.
We also define $F_5=F_2\circ G_1$ and $F_6=F_3\circ G_4$.
Then, the resulted nested fractal $K$ is called Pentakun and a subset $L$ is a Koch curve (see Figure~1).
The Hausdorff dimensions of $K$ and $L$ are $(\log5)/(\log \alpha)$ and $(\log4)/(\log \alpha)$,
respectively. There exists a canonical Dirichlet form $(\ce,\cf)$ on $L^2(K,\mu)$ 
where the corresponding diffusion is the fractional diffusion (see \cite{bar,kum0,lind} etc.),  
so $\cf=\Lambda_{2,\infty}^{d_w/2}
(K)$. It is known (see \cite{kum0}) that $d_w=(\log \frac{\sqrt {161}+9}{2})/(\log \alpha)$ 
and we can check all the assumptions similarly to the case of the Sierpinski gasket. 
Note that $C_2$ given in (B2)~(4) is equal to 4.
Thus, by Theorem~\ref{theo:mainthm1}, Theorem~\ref{theo:mainthm2} and Remark~\ref{rem:2.5},
\[
  \cf|_L=\Lambda_{2,2}^\beta (L)\qquad\mbox{where }~~\beta=\frac{d_w}{2}-\frac{\log 5-\log 4}{2\log \alpha}.
\]
For the Pentakun $K$, let $I'=\{2,3\}$. Then the corresponding self-similar subset $L'$ is a Cantor set
with Hausdorff dimension $(\log 2)/(\log \alpha)$. In this case we should set 
$\hat I=\{2,3,5,6\}$, so $I\ne \hat I$.
Again we can check all the assumptions similarly, so 
by Theorem~\ref{theo:mainthm1}, Theorem~\ref{theo:mainthm2} and Remark~\ref{rem:2.5},
\[
  \cf|_{L'}=\Lambda_{2,2}^{\beta'} (L')\qquad\mbox{where }~~\beta'=\frac{d_w}{2}-\frac{\log 5-\log 2}{2\log
\alpha}.
\]

In general, if $K$ is a nested fractal satisfying (A3), then 
there is a canonical Dirichlet form on $L^2(K,\mu)$ where the corresponding diffusion is 
the fractional diffusion (see \cite{bar,kum0,lind} etc.).
Let $L$ be a self-similar subset of $K$ given in the manner in the first part of Section~2, 
and satisfying (A2). 
In most cases, all the assumptions 
except (B4) can be checked similarly to the case of 
the Sierpinski gasket, so that we can use Theorem~\ref{theo:mainthm1} and Theorem~\ref{theo:mainthm2} to
characterize the trace space
if (B4)
holds. However, there are cases where (B4) does not hold -- see 4).
\bigskip

\noindent 3) Sierpinski carpets: 
Let $H_0=[0,1]^n$, 
$n\ge 2$,  
and let $l \in \bn$, $l\geq 2$ be fixed.
Set ${\cal Q}=\{\Pi_{i=1}^n [(k_i-1)/l,k_i/l]: 1\le k_i\le l,~k_i\in \bn~
(1\le i\le n)\}$, let $N \le l^n$ and $W=S=\{1,\ldots,N\}$.
Let $F_i$, $i\in S$
be orientation preserving affine maps of
$H_0$ onto some element of ${\cal Q}$. 
We assume that 
the sets $F_i(H_0)$ are distinct.
Set $H_1=\bigcup_{i\in I}F_i (H_0)$. 
Then the unique non-void compact set $K$ which satisfies $K=\bigcup_{i=1}^N F_i(K)$ 
is called the generalized Sierpinski carpet if the following holds: 
\begin{enumerate}
\item[(SC1)] (Symmetry) $H_1$ is preserved by all the isometries of the unit cube $H_0$.
\item[(SC2)] (Connected) $H_1$ is connected. 
\item[(SC3)] (Non-diagonality) Let $B$ be a cube in $H_0$ which is the union of
$2^n$ distinct elements of ${\cal Q}$. (So $B$ has side length $2l^{-1}$.) Then
if $\mbox{Int} (H_1 \cap B)$ is non-empty, it is connected.
\item[(SC4)] (Borders included) $H_1$ contains the line segment $\{x: 0\leq x_1
\leq 1, x_2=\cdots =x_n=0 \}$.
\end{enumerate}
Here (see \cite{bb}) (SC1) and (SC2) are essential, while
(SC3) and (SC4) are included for technical convenience.
The Sierpinski carpets are infinitely ramified: the critical
set $C_K$ in (\ref{eq:critpc}) is 
an infinite set,
and $K$ cannot be disconnected by removing a finite number of points.

 It is known (see \cite{bb,bbk,kz} etc.) 
that there is a self-similar Dirichlet form $(\ce,\cf)$ on $L^2(K,\mu)$ 
where the corresponding diffusion is the fractional diffusion.
In particular, $\cf=\Lambda_{2,\infty}^{d_w/2} (K)$ where $d_w=(\log \rho N)/(\log l)$,
$\rho$ given in (A5). 
Let $\fg=\{\mbox{the identity map}\}$ and $L=([0,1]^{n-1}\times\{0\})\cap K$ (cf. Figure 1).
Let $I=\{i\in S: F_i(K)\cap L\ne\emptyset\}$, $N_I=\# I$, and assume
\begin{equation}\label{karho}
\rho N_I>1.
\end{equation}
For simplicity, we assume that the $(n-1)$-dimensional Sierpinski carpet $L$ also satisfies the conditions corresponding to (SC1)--(SC4).
Then, (A1)--(A6) and (C1)--(C2) are easy to check with $M=1$. 
(A7) holds by (\ref{karho}), because 
$\beta=d_w/2-(d_f-d)/2=(\log \rho N_I)/(2\log l)$. 
(A8) holds by Proposition \ref{thm:pr4.2}. It is known that the corresponding self-adjoint
operator has compact resolvents (see \cite{bb,bbk,kz} etc.), so (B1) holds. 
Letting $\hat I=I$, we can check (B2).
For $w\in I^{m}$, $m\in\bz_+$,
let $x_0(w)\in[0,1]^n$ be the center of $K_w$ 
and $\Lambda_k(w)$ the intersection of $K$ and a cube in $\br^{n}$ with center $x_0(w)$ and length $(2k+1)l^{-m}$ for $k\in\bn$.
In order to assure (B3), 
assume for the moment 
that there exists some $k\ge6$ such that $\Lambda_k(w)$ is connected for all $w\in \bigcup_{m\in\bz_+}I^{m}$.
Let $l_0=(2k+1)n$ and take $m_0\in\bn$ such that $l^{m_0}\ge 2k+1$.
For each $w\in I^{m+m_0}$, $m\in\bz_+$, take $D'''(w)\subset D''(w)\subset D^\sharp(w)\subset D^{(1)}(w)\subset D^{(2)}(w)\subset D'(w)$ so that $K_{D'''(w)}=\Lambda_1(w)$, $K_{D''(w)}=\Lambda_2(w)$, $K_{D^\sharp(w)}=\Lambda_3(w)$, $K_{D^{(1)}(w)}=\Lambda_4(w)$, $K_{D^{(2)}(w)}=\Lambda_5(w)$, and $K_{D'(w)}=\Lambda_k(w)$.
With the use of Proposition~\ref{prop:ehi} and Proposition~\ref{prop:wv}, (B3) can be checked.
Here, the Harnack inequalities (EHI1) and (EHI2) are assured by \cite{bb,bbk,kz}. 
To be more precise, let  $\hat K=\bigcup_{x\in\{-1,0\}^n}(K+x)$, which is a subset of $[-1,1]^n$. 
Then, one can construct a Dirichlet form on $\hat K$ whose corresponding diffusion is the fractal diffusion 
in the same way as in \cite{bb,bbk,kz}. Indeed, $\hat K$ has enough symmetry for the coupling arguments in \cite{bb} 
to work. In this way, the Harnack inequalities (EHI1) and (EHI2) are assured. 
If for each $k$, there exists $w\in \bigcup_{m\in\bz_+}I^{m}$ such that $\Lambda_k(w)$ is not connected, then
take the connected component of $\Lambda_k(w)$ 
including $K_w$ in place of $\Lambda_k(w)$
and discuss similarly as above. By the covering argument, we
can check (B3).
(B4) is confirmed by an argument similar to the case of Sierpinski gaskets.
Thus, we have by Theorem~\ref{theo:mainthm1} and  Theorem~\ref{theo:mainthm2}, 
\[\cf|_L={\hat \Lambda}_{2,2}^\beta (L)\qquad\mbox{where }~~\beta=\frac{d_w}2-\frac 12 
\left(\frac{\log N}{\log l}-\mbox{dim}_H L\right).\]
Note that when $\partial [0,1]^n\subset K$, then $0<\beta<1$, so (\ref{karho}) holds
and $\cf|_L=\Lambda_{2,2}^\beta (L)$ by Remark \ref{rem:2.5}.  
Indeed, let $K_2=[0,1]^{n}$ and $K_1$ be a generalized Sierpinski carpet in $\br^{n}$  
with $\partial [0,1]^n\subset K_1$, 
which is 
determined by $\{F_i\}_{i}$
where $F_i([0,1]^{n})\cap \partial [0,1]^{n}\ne 
\emptyset$ for all $i$. Clearly, $K_1\subset K\subset K_2$. For each $K_i$, one can construct the 
self-similar Dirichlet form. Let $\rho_i$ be the scaling factor given in (A5).
By the shorting and cutting laws for electrical networks (see \cite{ds}), $\rho_2\le \rho\le 
\rho_1$. Then, $\rho_2=l^{2-n}$ and 
\begin{equation}\label{bbreses}
\frac{2}{l^{n-1}}+\frac {l-2}{l^{n-1}-(l-2)^{n-1}}\le \rho_1\le
\frac l{l^{n-1}-(l-2)^{n-1}},
\end{equation}
due to (5.9) in \cite{bb}. Since $L=[0,1]^{n-1}\times\{0\}$ and $N_I=l^{n-1}$ in this case,
we have $\rho N_I\ge \rho_2 N_I=l\ge 2$, so (\ref{karho}) holds and $\beta>0$.
Using (\ref{bbreses}),
\[\rho N_I\le \rho_1 N_I\le \frac {l^n}{l^{n-1}-(l-2)^{n-1}}<l^2,\]
where the last inequality is a simple computation. Thus $\beta<1$.

\medskip

\noindent 4) Vicsek sets: Let $a_1=(0,0), a_2=(1,0), a_3=(1,1), a_4=(0,1), a_5=(1/2,1/2)$
be points in $\br^2$ and define $F_i(x)=(x-a_i)/3+a_i$ for $x\in \br^2$ and $i=1,\ldots, 5$.
The unique non-void compact set $K$ which satisfies $K=\bigcup_{i=1}^5 F_i(K)$ 
is the Vicsek set. As in the case of 1), 
there is a self-similar Dirichlet form $(\ce,\cf)$ on $L^2(K,\mu)$ with $\rho=3$, 
where the corresponding diffusion is the fractional diffusion. In particular,
$\cf=\Lambda_{2,\infty}^{d_w/2} (K)$ where $d_w=(\log 15)/(\log 3)$.  
Let $L$ be the line segment from $(0,0)$ to $(1,1)$. Then, (B4) does not
hold. In this case, the trace of the Brownian motion on the Vicsek set is the Brownian motion
on the line segment.
Indeed, one can easily check condition $(H_1) - (H_3)$ in Section~8 of \cite{BP}
on the $1$-dimensional Sierpinski gasket which is the line. So, by  \cite{BP} Theorem 8.1,
one sees that the trace of the Brownian motion  on the Vicsek set 
is a constant time change of the Brownian motion on the line. 
Thus,   
\[\cf|_L=\Lambda_{2,\infty}^1 (L),\]
which is larger than $\Lambda_{2,2}^1 (L)$. This shows that (B4) is necessary for 
Theorem \ref{theo:mainthm1}. 
\section{Application: Brownian motion penetrating fractals}
In \cite{hk}, one of the authors constructed Brownian motions on {\it fractal fields}, a collection
of fractals with (in general) different Hausdorff dimensions (see also \cite{kum}). 
They are diffusion processes which 
behave as the appropriate fractal diffusions within each fractal component of the field and 
they penetrate each fractal. In \cite{hk}, 
a restricted assumption (Assumption 2.2 in \cite{hk}) was needed to construct such processes 
because we did not know the corresponding function spaces.
Our result in this paper can
be applied here and we can construct such penetrating diffusions  without the restricted assumption. 

Let $A_0$ be a countable set and 
let $\{K_i\}_{i\in A_0}\subset \br^n$ be a family of self-similar sets 
together with strong local, regular, and self-similar Dirichlet forms 
$(\ce_{K_i},\cf_{K_i})$ on $\bl^2(K_i,\mu_i)$, where $K_i$ and $\mu_i$ lie in the framework of Section~2.
We also regard $\mu_i$ as a measure on $\br^n$ by letting $\mu_i(\br^n\setminus K_i)=0$. 
We set $G=\bigcup_{i\in A_0} K_i$.

Let $A_1$ be another countable set and 
let $\{D_j\}_{j\in A_1}\subset \br^n$ be a family of disjoint domains
in $\br^n\setminus G$. 
Denote the closure of $D_j$ in $\br^n$ by $K_j$
and the Lebesgue measure restricted on $K_j$ by $\mu_j$.
Define $\tg=G\cup \bigl(\bigcup_{j\in A_1}K_j\bigr)$.
$\tg$ is called a {\it fractal field} generated by 
$\{K_i\}_{i\in A_0}$ and $\{D_j\}_{j\in A_1}$. 
(When $G$ is connected as in the introduction, we also call
$G$ a {\it fractal field} or a {\it fractal tiling}.) 

Denote by $A$ the disjoint union of $A_0$ and $A_1$. 
For $i,j\in A$ with $i\ne j$, let
$\Gamma_{ij}=K_i\cap K_j$.
Define $\Gamma=\bigcup_{i,j\in A,\,i\ne j}\Gamma_{ij}$. 
For $x\in \Gamma$, let $J_x:=\{i\in A: x\in K_i\}$ 
and define $N_x:=\bigcup_{i,j\in J_x,\,i\ne j}\Gamma_{ij}$.  
Throughout this section, we impose
the following assumption.\medskip 

\noindent{\bf Assumption A}
(1) For each compact set $C\subset \br^n$, 
$\#\{i\in A: C\cap K_i\ne \emptyset\}<\infty$.\\
(2) For each $i\in A_1$, $K_i\setminus D_i$ is a null set with respect to the Lebesgue measure on $\br^n$.
\medskip

For each $i\in A_1$, define 
${\cal D}(\ce_{K_i})=\{u\in C_0(K_i): u|_{D_i}\in W^{1,2}(D_i)\}$ and 
\[
\ce_{K_i} (u,v)=\frac 12\int_{D_i}(\nabla u(x),
\nabla v(x))_{\br^n}\,dx,~~~\mbox{for}~~u,v\in {\cal D}(\ce_{K_i}).
\]
Then, $(\ce_{K_i},{\cal D}(\ce_{K_i}))$ is closable on $L^2(K_i,\mu_i)$.
Its closure will be denoted by $(\ce_{K_i},\cf_{K_i})$.
It is easy to see that $(\ce_{K_i},\cf_{K_i})$ is a strong local regular Dirichlet form.

For $x\in \Gamma$ and $i\in J_x$,
define $\beta_{x,i}=d_w(K_i)/2-(d_f(K_i)-d_f(N_x\cap K_i))/2$.
Here, $d_w(K_i)$ is defined in (A7) for $(\ce_{K_i},\cf_{K_i})$ if $i\in A_0$ and $d_w(K_i)$ is defined as $2$ if $i\in A_1$,
and $d_f(K_i)$ and $d_f(N_x\cap K_i)$ are the Hausdorff dimensions of $K_i$ and $N_x\cap K_i$, respectively.

We will also assume the following throughout this section. 

\medskip 
\noindent {\bf Assumption B }
(1) For $i\in A_0$, 
$(\ce_{K_i},\cf_{K_i})$ is a 
strong local regular Dirichlet form on 
$\bl^2(K_i,\mu_i)$ which satisfies (A1),
(A3), the first identity of (A4),
(A5), and (A6) in Section~2.\\
(2) For each $x\in\Gamma$ and $i\in J_x\cap A_0$, $N_x\cap K_i$ is a finite number of union of compact self-similar sets $\{L_j\}$
that are constructed by the same number of contraction maps and each of which satisfies (A2), the second identity of (A4),
(A7), and (A8) in Section~2. 
Further, (C1)$^*$--(C2)$^*$ in Remark~\ref{thm:finunil}
holds with $K=K_i$ and $L=N_x\cap K_i$.\\
(3) For each $x\in \Gamma$ and $i\in J_x\cap A_1$, $N_x\cap K_i$ is a closed Alfors $d_{x,i}$-regular set with some $d_{x,i}$.\\
(4) 
For every $x\in \Gamma$, $\beta_{x,i}>0$ for all $i\in J_x$, and the set $\Lambda_x:=\{f\in  C_0(N_x): f|_{N_x\cap K_i}\in \Lambda_{2,2}^{\beta_{x,i}}(N_x\cap K_i)\mbox{ for all }i\in J_x\}$ is dense in $C_0(N_x)$.
\medskip

We will give several remarks.
When $i\in A_1$, we have $d_f(K_i)=n$ and $d_f(K_i\cap N_x)=d_{x,i}$.
The set $\Lambda_x$ is closed under the operation of the normal contraction; $0\vee f\wedge1\in \Lambda_x$ for $f\in \Lambda_x$.
If $N_x$ itself is an Alfors regular set and $\beta_{x,i}\in (0,1)$ for all $i\in J_x$, 
then $\Lambda^{\max_{i\in J_x}\beta_{x,i}}_{2,2}(N_x)\cap C_0(N_x)$ (which is a subset of $\Lambda_x$) is dense 
in $C_0(N_x)$ by Chapter~V, Proposition~1 in \cite{jw} and Theorem~3 in \cite{stos}.
The condition $\beta_{x,i}\in (0,1)$ holds, for example, if $i\in N_x\cap A_1$ and $d_{x,i}\in(n-2,n)$, because then $\beta_{x,i}=1-(n-d_{x,i})/2\in (0,1)$.

Define a measure $\tmu$ on $\tg$ by $\tmu=\sum_{i\in A}\mu_i$.
We now define a bilinear form $(\tce,{\cal D}
(\tce))$ on $\bl^2 (\tg,\tmu)$ as follows:
\begin{eqnarray*}
\tce (u,v)&=& \sum_{i\in A} \ce_{K_i} (u|_{K_i},v|_{K_i})
~~\mbox{for}~u,v \in {\cal D} (\tce),\\ 
{\cal D} (\tce) &=& \{u\in C_0(\tg): u|_{K_i}\in\cf_{K_i}
\mbox{ for all }i\in A \mbox{ and } \tce (u,u)<\infty\}.
\end{eqnarray*}
Then, the following is easy to check.
\begin{lem}
\begin{enumerate}
\item[$(1)$] $(\tce,{\cal D} (\tce))$ is closable in $\bl^2 (\tg,\tmu)$.
\item[$(2)$] ${\cal D} (\tce)$ is an algebra.
\item[$(3)$] For $i\in A$, $x\in K_i$, and for $U(x)$ which is a 
neighborhood of $x$ in $K_i$, there exists $f\in \cf_{K_i}\cap C_0(K_i)$ 
such that $f(x)>0$ and
$\supp f\subset U(x)\cap K_i$, where $\supp f$ denotes the
support of $f$.
\end{enumerate}
\label{theo:easy}\end{lem}
Now, let $(\tce, \tcf)$ be the closure of $(\tce, {\cal D} (\tce))$. 
We then have the following.

\begin{theo}\label{theo:maincon}
$(\tce,\tcf)$ is a strong local 
regular Dirichlet form on $\bl^2 (\tg,\tmu)$.
\end{theo}

Note that the strong local property of $(\tce,\tcf)$ can be 
easily deduced from those of the original forms on 
$\{K_i\}_{i\in A}$. Therefore, it is enough to prove the
regularity of $(\tce,\tcf)$. 
For the proof of it, 
the key part is to prove the following.
\begin{pro}\label{theo:keypro}
\begin{enumerate}
\item[$(1)$] For each $x\ne y\in \tg$, there exists $g\in {\cal D}(\tce)$
such that $g(x)\ne g(y)$.
\item[$(2)$] For any compact set $L$ in $\tg$, there exists $f\in {\cal D}(\tce)$
such that $f=1$ on $L$.
\end{enumerate}
\end{pro}
Once this proposition is established, it is easy to prove 
the regularity of $(\tce,\tcf)$
(see \cite{hk}), so we will only prove the proposition. 
\Proofwithname{Proof of Proposition \ref{theo:keypro}.}
Let $B(x,r)$ denote the open ball in $\br^n$ with center $x\in\br^n$ and radius $r$. 
When either $x$ or $y$ is in the compliment of $\Gamma$,
then (1) is clear by Lemma \ref{theo:easy}~(3), so we will consider the case 
$x,y\in \Gamma$. By Assumption~A~(1), $\# J_x<\infty$.
Since each $K_j$ is closed, by Assumption~A~(1), there exists 
$r_x>0$ such that $B(x,r_x)\cap K_j\ne \emptyset$ 
if and only if $j\in J_x$, and $y\notin B(x,r_x)$. 
Since $\Lambda_x$ is dense in $C_0(N_x)$ by Assumption~B~(4),
there exists $u\in \Lambda_x$ such that 
$u|_{B(x,r_x/2)}=1$ and $u|_{B(x,3r_x/4)^c}=0$.

Now, by Assumption~B~(1), (2) and the extension theorem (Remark~\ref{thm:finunil}), 
for each $i\in J_x\cap A_0$, 
there exists ${\hat u}_i\in \cf_{K_i}\cap C(K_i)$ such that 
${\hat u}_i|_{N_x\cap K_i}=u$. 
For each $i\in J_x\cap A_1$, since $N_x\cap K_i$
is a closed Alfors $d_{x,i}$-regular set, we have 
\begin{equation}\label{dsetra}
W^{1,2}(\br^n)|_{N_x\cap K_i}=\Lambda_{2,2}^{1-(n-d_{x,i})/2}(N_x\cap K_i)
\end{equation}
(see \cite{jw}). By carefully tracing the proof of the extension 
theorem in (\ref{dsetra}), we see that there exists 
${\hat u_i}\in W^{1,2}(\br^n)\cap C_0(\br^n)$ such that 
${\hat u_i}|_{N_x\cap K_i}=u$  
(see, for instance, pages 77--78 in \cite{kum}).
For both cases, since $(\ce_{K_i},\cf_{K_i})$ is regular,  
by multiplying a function in $\cf_{K_i}\cap C_0(K_i)$ which is
$1$ in $B(x,3r_x/4)$ and $0$ outside $B(x,r_x)$, 
we may assume $\supp {\hat u_i}\subset B(x,r_x)$.
Define $g\in C_0(\tg)$ as $g|_{K_i}={\hat u}_i$ for $i\in J_x$ and 
$g|_{K_i}\equiv 0$ otherwise. Then, $g\in {\cal D} (\tce)$, 
$g(x)=1$ and $g(y)=0$. We thus obtain the desired function.

The proof of (2) is quite similar, so we omit it 
(see Proposition 2.6~(2) in \cite{hk}).
\QED
Denote the 1-capacity associated with $(\ce_{K_i},\cf_{K_i})$ and $(\tce,\tcf)$ by $\Cp_{K_i}$ and $\Cp_{\tilde G}$, respectively.
By definition, it is easy to see that $u|_{K_i}\in \cf_{K_i}$ for any $i\in A$ and $u\in\tcf$.
Further, $\Cp_{K_i}(H)\le \Cp_{\tilde G}(H)$ for any $i\in A$ and $H\subset K_i$. 
For $i\in A$, let $\cf_{K'_i}=\bigl\{f\in \cf_{K_i}:\tilde f=0 \mbox{ q.e.\ on } \Gamma\bigr\}$ 
and $\tcf_{i}=\bigl\{f\in \tcf:\tilde f=0 \mbox{ q.e.\ on }\bigcup_{j\in A\setminus\{i\}}K_j\bigr\}$, where $\tilde f$ is a (corresponding) quasi-continuous modification of $f$.

We will denote by $(\{\tx_t\}_{t\ge 0},\{\tilde P_x\}_{x\in\tilde G})$ the diffusion process corresponding to $(\tce,\tcf)$. 
The following proposition shows that $\{\tilde X_t\}$ behaves on $K_i$
in the same way as the diffusion process associated with 
$(\ce_{K_i},\cf_{K_i})$ until the process hits $\Gamma$.
\begin{pro}\label{theo:part}
$(\ce_{K_i},\cf_{K'_i})$ and $(\tce,\tcf_i)$ give the same Dirichlet forms 
on $L^2(K_i,\mu_i|_{K_i\setminus\Gamma})$,  
by identifying the measure space $(\tilde G,\mu_i|_{K_i\setminus\Gamma})$ with 
$(K_i,\mu_i|_{K_i\setminus\Gamma})$. 
In particular, the corresponding parts of the processes on $K_i\setminus\Gamma$ are the same.
\end{pro}
\Proof
 It is easy to see that $f\in\tcf_i$ satisfies that $f|_{K_i}\in \cf_{K'_i}$, so we will prove the converse.
 Let $f\in \cf_{K'_i}$.
 By Theorem~4.4.3 of \cite{fot}, we can take an approximation sequence of $f$ from $\cf_{K'_i}\cap C_0(K_i\setminus\Gamma)$. 
 Therefore, the $0$-extension of $f$ outside $K_i$ is an element of $\tcf_i$. 
\QED
For each distinct $i,j\in A$, we denote $K_i\sim K_j$ 
if $\Cp_{K_l}(\Gamma_{ij})>0$ for $l=i$ and $j$. 
We now assume the following in addition to 
Assumptions~A and B. 
\medskip

\noindent{\bf Assumption C}
(1) For each $i\in A_0$, 
$(\ce_{K_i},\cf_{K_i})$ is irreducible. \\
(2) For each distinct $i,j\in A$, there exist $k\in\bn$ and
a sequence $i_0,i_1,\ldots, i_k\in A$ such that 
$K_{i_0}=K_i$, $K_{i_k}=K_j$ and $K_{i_l}\sim K_{i_{l+1}}$ 
for $l=0,1,\ldots, k-1$.
\\
(3) For each distinct $i,j\in A$ with $K_i\sim K_j$, there exists a positive Radon measure $\nu_{ij}$ on $\Gamma_{ij}$ such that $\nu_{ij}(\Gamma_{ij})>0$ and $\nu_{ij}$ is smooth with respect to both $(\ce_{K_i},\cf_{K_i})$ and $(\ce_{K_j},\cf_{K_j})$.
\medskip
 
Note that, when $i\in A_1$, $(\ce_{K_i},\cf_{K_i})$ is irreducible since $D_i$ is connected.
(See, e.g.\ Theorem~4.5 in \cite{ou}, for the proof.)

For each nearly Borel set $B\subset \br^n$, define
$\sigma_B=\inf \{t>  0: \tilde {X}_t\in B\}$. 
The next proposition shows that ${\tilde X}_t$ penetrates into each $K_i$. 
\begin{pro} The following holds  
for any nearly Borel set $B$ with 
$\Cp_{\tg}(B)>0$.
\begin{equation}
\tilde {P}^x (\sigma_B<\infty)>0
\mbox { for $(\tce,\tcf)$-quasi every }x\in \tilde G.
\label{eq:hitinqq}\end{equation} 
Especially, if $B$ is 
a subset of a certain $K_i$ 
with $\Cp_{K_i}(B)>0$, then
\rom{(\ref{eq:hitinqq})} holds.
\end{pro}
\Proof 
By virtue of Theorem 4.6.6 in \cite{fot}, it is enough to prove that $(\tce,\tcf)$ is irreducible. 
We first recall the following fact.
Let $(\ce,\cf)$ be a local Dirichlet form. (Here, the locality means $\ce(f,g)=0$ if $fg=0$ a.e. All Dirichlet forms appearing in this article are local in this sense; see \cite{sch}.)
Let $Y$ be a measurable subset of the state space and  $\cc$ a dense set in $\cf$.
Then,  $Y$ is an invariant set if and only if $1_Y\cdot u\in\cf$ for any $u\in\cc$.
This is verified by Theorem~1.6.1 in \cite{fot} and a usual approximation argument. 

Now, let $M$ be an invariant set for $(\tce,\tcf)$.
Fix $i\in A$ and take $u\in \cf_{K_i}\cap C_0(K_i)$.
We can take $v\in{\cal D}(\tce)$ such that $v=1$ on $\supp u$ by Proposition~\ref{theo:keypro}~(2).
Then, $1_M\cdot v\in\tcf$, which implies that $(1_{M}\cdot v)|_{K_i}\in \cf_{K_i}$.
Therefore, $u\cdot(1_{M}\cdot v)|_{K_i}=u\cdot1_{M\cap K_i}$ also belongs to $\cf_{K_i}$.
Since $\cf_{K_i}\cap C_0(K_i)$ is dense in $\cf_{K_i}$, we obtain that $M\cap K_i$ is an invariant set for $(\ce_{K_i},\cf_{K_i})$.
By the irreducibility of $(\ce_{K_i},\cf_{K_i})$, either $\mu_i(M\cap K_i)=0$ or $\mu_i(K_i\setminus M)=0$ holds.

By this argument, there exists a subset $A'$ of $A$ such that $M=\bigcup_{i\in A'}K_i$ $\tmu$-a.e.
Assume that $M$ is a nontrivial invariance set.
Then, $A'\ne\emptyset$, $A'\ne A$, and there exist $i\in A'$ and $j\in A\setminus A'$ such that $K_i\sim K_j$ by Assumption~C~(2).
Take a compact set $H\subset \Gamma_{ij}$ such that $\nu_{ij}(H)>0$, and a relatively compact open set $H'$ including $H$.
Take $v\in{\cal D}(\tce)$ such that $v=1$ on $H'$
and let $u=1_M\cdot v\in\tcf$.
Denote by $\tilde u$ the quasi-continuous modification of $u$ w.r.t.\ $(\tce,\tcf)$.
Then, $\tilde u|_{K_l}$ is also quasi-continuous w.r.t.\ $(\ce_{K_l},\cf_{K_l})$ for $l=i,j$.
Since $\tilde u=1$ $\mu$-a.e.\ on $H'\cap K_i$, we have $\tilde u=1$ $\ce_{K_i}$-q.e.\ on $H\subset H'\cap K_i$. 
By Assumption~C~(3), $\tilde u=1$ $\nu_{ij}$-a.e.\ on $H$. 
On the other hand, since $\tilde u=0$ $\mu$-a.e.\ on $H'\cap K_j$, we have $\tilde u=0$ $\ce_{K_j}$-q.e.\ on $H$. 
Therefore, $\tilde u=0$ $\nu_{ij}$-a.e.\ on $H$. 
This is a contradiction, which deduces that $(\tce,\tcf)$ is irreducible.
\QED
The fractal field in Figure 2 satisfies 
Assumptions~A, B and C, so there is a penetrating 
diffusion on the field.

In \cite{hk}, detailed properties of $\tilde {X}_t$ such as heat kernel bounds 
and large deviation estimates are established under strong assumptions 
such as Assumption $2.2$ in \cite{hk}. Using the results given in this section,
one can relax the assumption and obtain the same results by the same proof given 
in \cite{hk}, 
when each Dirichlet form is the resistance form in the sense of \cite{kig}. 
\bigskip

\noindent
{\sc Acknowledgements.} The authors thank Professor R.~S. Strichartz 
for fruitful discussions and valuable comments.


\begin{thebibliography}{10}
\bibitem{adh}
D.~R.~Adams and L.~I.~Hedberg,
\newblock Function spaces and potential theory,
\newblock Springer, Berlin Heidelberg, 1996. 

\bibitem{bar}
M.~T.~Barlow,
\newblock Diffusions on fractals,
\newblock Lectures in Probability Theory and Statistics:
Ecole d'\'et\'e de probabilit\'es de Saint-Flour XXV, Springer, New York, 1998.

\bibitem{bb}
M.~T.~Barlow and R.~F.~Bass,
\newblock Brownian motion and harmonic analysis on Sierpinski carpets,
\newblock {\it Canad. J. Math.}, {\bf 51} (1999), 673--744.

\bibitem{bbk}
M.~T.~Barlow, R.~F.~Bass and T.~Kumagai, 
\newblock Stability of parabolic Harnack inequalities on metric measure spaces, 
\newblock {\it J. Math. Soc. Japan}, to appear.

\bibitem{bk}
M.~T. Barlow and T. Kumagai,
\newblock Transition density asymptotics for some diffusion processes 
with multi-fractal structures,
\newblock {\it Electron. J. Probab., \bf 6} (2001), 1--23.

\bibitem{BP} 
M. T. Barlow and E. A. Perkins,
\newblock Brownian Motion on the Sierpinski Gasket,
\newblock {\it Probab. Theory Relat. Fields, \bf 79} (1988), 543--623.

\bibitem{ben} 
A. Ben Amor, 
\newblock Trace inequalities for operators associated to regular Dirichlet forms,
\newblock {\it  Forum Math., \bf 16}  (2004), 417--429.

\bibitem{bod} 
M. Bodin, 
\newblock Characterisations of function spaces on fractals,
\newblock Ume\aa~University, Department of Math. and Math. Stat.,
Doctoral Thesis, No. 32, 2005

\bibitem{ds} 
P. G.~Doyle and J. L.~Snell,
\newblock Random walks and electrical networks,
\newblock Washington, Math. Assoc. of America, 1984.

\bibitem{fot}
M. Fukushima, Y. Oshima, and M. Takeda,
\newblock Dirichlet forms and symmetric Markov processes, 
\newblock de Gruyter, Berlin, 1994.

\bibitem{gri}
A. Grigoryan,
\newblock Heat kernels and function theory on metric measure spaces,
\newblock {\it Contemp. Math., \bf 338} (2003), 143--172.

\bibitem{hk}
B. M. Hambly and T. Kumagai,
\newblock Diffusion processes on fractal fields: heat kernel estimates and large 
deviations, 
\newblock {\it Probab. Theory Relat. Fields, \bf 127} (2003), 305--352.

\bibitem{hin}
M. Hino, 
\newblock On singularity of energy measures on self-similar sets, 
\newblock {\it Probab. Theory Relat. Fields, \bf 132} (2005), 265--290.

\bibitem{j2}
A. Jonsson,
\newblock Brownian motion on fractals and function spaces, 
\newblock {\it Math. Z., \bf 222} (1996), 496--504.

\bibitem{jo}
A.~Jonsson,
\newblock A trace theorem for the Dirichlet form on the Sierpinski gasket,
\newblock {\it Math. Z.}, to appear. 

\bibitem{jw}
A.~Jonsson and H.~Wallin, 
\newblock Function spaces on subsets of $\br^n$,
\newblock Mathematical Reports, Vol.~{\bf 2}, Part $1$ (1984), Acad. Publ., Harwood.

\bibitem{kam}
A. Kamont,
\newblock A discrete characterization of Besov spaces,
\newblock {\it Approx. Th. Appl., \bf 13} (1997), 63--77. 

\bibitem{kig}
J.~Kigami,
\newblock Analysis on fractals,
\newblock Cambridge Univ. Press, Cambridge, 2001.

\bibitem{kum0}
T. Kumagai,
\newblock Estimates of transition densities for Brownian motion on nested
fractals,
\newblock {\it Probab. Theory Relat. Fields, \bf 96} (1993), 205--224.

\bibitem{kum}
T. Kumagai,
\newblock Brownian motion penetrating fractals
-An application of the trace theorem of Besov spaces-,
\newblock {\it J. Func. Anal., \bf 170} (2000), 69--92.

\bibitem{kz}
S.~Kusuoka and X. Y.~Zhou,
\newblock Dirichlet forms on fractals: Poincar\'e constant and resistance,
\newblock {\it Probab. Theory Relat. Fields, \bf 93} (1992), 169--196.

\bibitem{lei}
L.~Leindler,
\newblock Generalization of inequalities of Hardy and Littlewood,
\newblock{\it Acta Sci. Math. (Szeged), \bf 31} (1970), 279--285.

\bibitem{lind}
T.~Lindstr\o m,
\newblock Brownian motion on nested fractals,
\newblock Memoirs Amer. Math. Soc. $420$, {\bf 83}, 1990.

\bibitem{ou}
E. M.~Ouhabaz,
\newblock Analysis of heat equations on domains,
\newblock London Mathematical Society Monographs Series, 31. Princeton University Press, Princeton, NJ, 2005.

\bibitem{sch}
B. Schmuland, 
\newblock On the local property for positivity preserving coercive forms,
\newblock  Dirichlet forms and stochastic processes (Beijing, 1993),  
345--354, de Gruyter, Berlin, 1995.

\bibitem{ste}
E. M. Stein,
\newblock Singular integrals and differentiability properties of functions, 
\newblock Princeton University Press, Princeton, 1970.

\bibitem{stos}
A. St\'os,
\newblock Symmetric $\alpha$-stable processes on $d$-sets,
\newblock {\it Bull. Polish Acad. Sci. Math., \bf 48} (2000), 
237--245.

\bibitem{stri}
R. S. Strichartz,
\newblock Function spaces on fractals, 
\newblock {\it J. Funct. Anal., \bf 198} (2003), 43--83. 

\bibitem{tr}
H. Triebel,
\newblock Fractals and Spectra, 
\newblock  Monographs in Math., Vol. 91,
Birkh\"auser, Basel-Boston-Berlin, 1997.

\end{thebibliography}
\end{document}